\documentclass[10pt, conference, compsocconf]{IEEEtran}
\ifCLASSINFOpdf
\else
\fi

\usepackage{graphicx}
\usepackage{subfigure} 
\usepackage{amsmath}
\usepackage{amsthm}
\usepackage{amssymb}
\usepackage{algorithm}
\usepackage{algorithmic}
\usepackage{url}
\usepackage{float}
\usepackage{multirow}
\usepackage{tikz}
\usepackage[numbers,sort&compress]{natbib}
\newtheorem{theorem}{Theorem}[section]
\newtheorem{lemma}[theorem]{Lemma}
\newtheorem{prop}[theorem]{Proposition}
\newtheorem{corollary}[theorem]{Corollary}

\def \R{\mathbb R}
\def\tU{\tilde{U}}
\def\la{\leftarrow}
\def\tV{\tilde{V}}
\def\tS{\tilde{\Sigma}}
\def\S{\Sigma}
\def\l{\lambda}
\def\s{\sigma}

\def\Diag{\text{Diag}}

\def\etal{et al.\;}

\newcommand{\NN}[1]{\| #1 \|_{*} }                   % nuclear norm
\newcommand{\FN}[1]{\| #1 \|_F }                     % frobenious norm
\newcommand{\SN}[1]{\| #1 \|_2 }                     % square norm

\newcommand{\Prox}[2]{\text{prox}_{#1}(#2)}
\newcommand{\SO}[1]{\mathcal{P}_{\Omega}(#1)}
\newcommand{\Tr}[1]{\text{tr}(#1)}
\newcommand{\Span}[1]{\text{span}(#1)}

\newcommand{\f}[1]{f(#1)}

\begin{document}
\title{Fast Low-Rank Matrix Learning with Nonconvex Regularization}


\author{\IEEEauthorblockN{Quanming Yao \quad James T. Kwok \quad Wenliang Zhong}
\IEEEauthorblockA{Department of Computer Science and Engineering \\
Hong Kong University of Science and Technology \\
Hong Kong\\
\{qyaoaa, jamesk, wzhong\}@cse.ust.hk}
%\and
%\IEEEauthorblockN{Leon Wenliang Zhong}
%\IEEEauthorblockA{line 1 (of Affiliation): dept. name of organization\\
%line 2: name of organization, acronyms acceptable\\
%line 3: City, Country\\
%line 4: Email: name@xyz.com}
}


%Telephone: (800) 555--1212, Fax: (888) 555--1212}


\maketitle

\begin{abstract}
Low-rank modeling has a lot of important applications in machine learning, computer vision and
social network analysis. While the matrix rank is often approximated by the convex nuclear
norm, the use of nonconvex low-rank regularizers has demonstrated better recovery performance.
However, the resultant optimization problem is
much more challenging.
A very recent state-of-the-art is based on the proximal gradient algorithm.
However, it requires an expensive full SVD in each proximal step.
In this paper, we show that 
for many commonly-used nonconvex low-rank regularizers, 
a cutoff can be derived to automatically threshold the singular
values obtained from the proximal operator. This allows the use of
power method to approximate the SVD efficiently.
Besides, the proximal operator can be reduced to that of a much smaller matrix
projected onto this leading subspace.
Convergence, with a rate of $O(1/T)$ where $T$ is the number of iterations,
can be guaranteed. Extensive experiments are performed on matrix completion and robust principal component analysis. The proposed method achieves
significant speedup over the state-of-the-art. Moreover, the matrix solution obtained is more
accurate and has a lower rank than that of the traditional nuclear norm regularizer.
\footnote{This is the long version of conference paper appeared in ICDM 2015 \cite{yao2015fast};
	Code is available at: \url{https://github.com/quanmingyao/FaNCL}.}
\end{abstract}


% \fi
\IEEEpeerreviewmaketitle

\section{Introduction}
\label{sec:intro}

The learning of 
low-rank matrices is a central issue in many machine learning problems. For example, 
matrix completion
\cite{candes2009exact},
which is one of the most successful approaches in collaborative filtering,
assumes that the target ratings matrix is low-rank. 
Besides 
collaborative filtering,
matrix completion has also been used on  tasks such as
sensor networks \cite{biswas-06},
social network analysis \cite{kim-11}, and image processing \cite{liu2013tensor,qyao2015color}.

Another important use of low-rank matrix learning is robust principal component analysis (RPCA) \cite{candes2011robust}, 
which assumes the target matrix is low-rank and also corrupted by  
sparse 
data noise.
It is now popularly used in various computer vision applications, such as shadow removal of aligned faces
and background modeling of surveillance videos \cite{candes2011robust,sun2013robust}.
Besides,
low-rank minimization has also been used in tasks such as
multilabel learning \cite{zhu2010image}
and multitask learning \cite{gong2012robust}.

However, rank minimization is NP-hard. To alleviate this problem, a common
approach is to use  instead
a convex surrogate such as 
the nuclear norm (which is the sum of singular values of the matrix).
It is known that the
nuclear norm is the tightest convex lower bound of the rank.
Besides, 
there are theoretical guarantees that 
the incomplete matrix 
can be recovered
with
nuclear norm regularization
\cite{candes2009exact,candes2011robust}. 
Moreover, though the nuclear norm is non-smooth, the resultant optimization problem can often be 
solved efficiently using modern tools such as
accelerated proximal gradient descent 
\cite{ji2009accelerated},
Soft-Impute \cite{mazumder2010spectral},
and active subspace selection methods \cite{hsieh2014nuclear}.

Despite its success, recently there have been numerous attempts that use
nonconvex surrogates to better
approximate the rank function.
The key idea is 
that the larger, and thus more informative, singular values
should be less penalized.
Example nonconvex low-rank regularizers include
the capped-$\ell_1$ penalty \cite{zhang2010analysis},
log-sum penalty (LSP) \cite{candes2008enhancing},
truncated nuclear norm (TNN) \cite{hu2013fast},
smoothly clipped absolute deviation (SCAD) \cite{fan2001variable}, and
minimax concave penalty (MCP) \cite{zhang2010nearly}. Empirically, these nonconvex regularizers achieve
better recovery performance than the convex nuclear norm regularizer.

However, the resultant nonconvex optimization problem is much more challenging.
One approach is to use the 
concave-convex procedure \cite{yuille-03}, which
decomposes the nonconvex regularizer
into a difference of convex functions 
\cite{zhang2010analysis,hu2013fast}.
However, a sequence of relaxed problems have to be solved,
and can be computationally expensive \cite{gong2013general}.
A more efficient method is the recently proposed iteratively re-weighted nuclear norm (IRNN) algorithm
\cite{canyi2014}. It is based on the observation that existing nonconvex regularizers are all
concave and their super-gradients are non-increasing. Though IRNN still has to iterate, each
of its iterations only involves computing the 
super-gradient of the regularizer and a singular value decomposition (SVD).
However, performing SVD on a $m\times n$ matrix (where $m\geq n$) 
still takes $O(mn^2)$
time, and can be expensive when the matrix is large.

Recently, the proximal gradient algorithm has 
also been used on this nonconvex low-rank minimization problem \cite{sun2013robust,hu2013fast,canyi2014,lu2015generalized}.
However, it requires performing the full SVD in each proximal step, 
which is expensive for large-scale applications. 
To alleviate this problem, we first observe that for the commonly-used nonconvex low-rank
regularizers,
the singular values obtained from the corresponding proximal operator can be automatically  
thresholded.
One then only needs to find the leading singular values/vectors in order to generate the next iterate. 
By using 
the power method \cite{halko2011finding},
a fast and accurate 
approximation of such a subspace can be obtained.
Moreover, instead of computing the proximal operator on a large matrix, one only needs to
compute that on its projection onto this leading subspace. The size of the matrix is  
significantly reduced and the proximal operator can be made much more efficient.
In the context of matrix completion problems,
further speedup is possible by
exploiting a special ``sparse plus low-rank'' structure of the matrix iterate.

The rest of the paper is organized as follows. 
Section~\ref{sec:review} reviews the related work.
The proposed algorithm
%, with analysis and extensions, 
is presented in Section~\ref{sec:proalg};
%Extensive 
Experimental results on matrix completion and RPCA
%robust principal component analysis 
are shown in 
Section~\ref{sec:expts}, and  the last section gives some
concluding remarks.

%\noindent
%Notation:
In the sequel, the transpose of vector/matrix is denoted by the superscript $(\cdot)^\top$. 
For a 
$m\times n$ matrix
$X$,
$\Tr{X}$ is its trace,
$\|X\|_F = \Tr{X^{\top}X}$ is the Frobenius norm, 
and $\NN{X} = \sum_{i=1}^m \sigma_i$ is the nuclear norm.
Given $x=[x_i] \in \R^{m}$, 
$\Diag(x)$ constructs a $m \times m$ diagonal matrix whose $i$th diagonal element is $x_i$.  
Moreover, $I$ denotes the identity matrix.
For a differentiable function $f$, we use $\nabla f$ for its gradient.
For a nonsmooth function,
we use $\partial f$ for its subdifferential.

%-----------------------------------------------------------------------------------------
% Review
%-----------------------------------------------------------------------------------------

\section{Background}
\label{sec:review}

%%%%%%%%%%%%%%%%%%%%%%%%%%%%%

\subsection{Proximal Gradient Algorithms}

In this paper, 
we consider composite optimization problems  of the form
\begin{equation} \label{eq:problem}
\min_x F(x) \equiv f(x) + \lambda r(x),
\end{equation} 
where $f$ is smooth and $r$ is nonsmooth.  
In many machine learning problems, $f$ is the loss and $r$ a low-rank regularizer.
In particular, we make the following assumptions on $f$.
\begin{itemize}
	\item[A1.] $f$, not necessarily convex, is differentiable with $\rho$-Lipschitz continuous
	gradient, i.e., $\FN{\nabla \f{X_1}-\nabla \f{X_2}} \le \rho \FN{X_1-X_2}$.
	Without loss of generality, we assume that $\rho \le 1$.
	
	\item[A2.] $f$ is bounded below, i.e., $\inf \f{X} > -\infty$.
\end{itemize}
In recent years, proximal gradient algorithms 
\cite{parikh2014proximal} have been widely used for solving
(\ref{eq:problem}).
At each iteration $t$, a quadratic function
is used to upper-bound the smooth $f$ at the current iterate $x^t$, while
leaving the nonsmooth $r$ intact. 
For a given stepsize $\tau$,
the next iterate $x^{t+1}$ is obtained as
\begin{eqnarray*}
	\lefteqn{\arg\min_x \nabla \f{x^t}^\top (x - x^t) + \frac{\tau}{2} \|x - x^t\|^2 + \l r(x)} \\
	& = & 
	\arg\min_x \frac{1}{2} \|x - z^t\|^2  + \frac{\l}{\tau} r(x)
	\equiv \Prox{\frac{\l}{\tau} r}{z^t},
\end{eqnarray*}
where $z^t = x^t - \frac{1}{\tau}\nabla \f{x^t}$, and
$\Prox{\frac{\l}{\tau} r}{\cdot}$ is the proximal operator \cite{parikh2014proximal}.
Proximal gradient algorithms 
can be further accelerated,
by replacing $z^t$ with a proper linear combination of $x^t$ and $x^{t-1}$.
In the sequel, as our focus is on learning low-rank matrices, $x$ in (\ref{eq:problem})
becomes a $m\times n$ matrix $X$.\footnote{In the following, we assume $m \le n$.}

%%%%%%%%%%%%%%%%%%%%%%%%%%%%%

\subsection{Convex and Nonconvex Low-Rank Regularizers}

An important factor for the success of proximal gradient algorithms is that its proximal operator 
$\Prox{\mu r}{\cdot}$ 
can be efficiently computed.
For example, 
for the 
nuclear norm $\NN{X}$,
the following Proposition shows that its proximal operator  
has a closed-form solution.

\begin{prop} \label{lem:svt} \cite{cai2010singular}
	$\Prox{\mu \NN{\cdot}}{X} 
	= U \left( \Sigma - \mu I\right)_{+} V^{\top}$,
	where $U \Sigma V^{\top}$ is the SVD of $X$, and $(Z)_{+} = [\max (Z_{ij}, 0)]$. 
\end{prop}

While the convex nuclear norm makes the low-rank 
optimization problem easier,
it may not be a good approximation of the matrix rank 
\cite{sun2013robust,hu2013fast,canyi2014,lu2015generalized}. As mentioned in
Section~\ref{sec:intro}, a number of nonconvex surrogates for the rank have been recently
proposed.  In this paper, we make the following assumption on the low-rank regularizer
$r$ in (\ref{eq:problem}).

\begin{itemize}
	\item[A3.] $r$ is possibly non-smooth and nonconvex, and 
	of the form 
	$r(X)=\sum_{i = 1}^m \hat{r}(\sigma_i)$,
	where $\sigma_1 \ge \dots \ge \sigma_m \ge 0$ are singular values of $X$,
	and $\hat{r}(\s)$ is a concave and non-decreasing function of $\s \ge 0$ with $\hat{r}(0) = 0$. 
\end{itemize}
All nonconvex low-rank regularizers introduced in Section~\ref{sec:intro}
satisfy this assumption.
Their corresponding $\hat{r}$'s are 
shown in 
Table~\ref{tab:lwregu}.

\begin{table}[H]
	\centering
	\renewcommand{\arraystretch}{1.2}
	\caption{$\hat{r}$'s for some popular nonconvex low-rank regularizers. 
		For the TNN regularizer, $\theta \in \{1,\dots,n\}$ is the number of leading singular values that are not penalized.}
	\vspace{-5px}
	\begin{tabular}{c | c} \hline
		& $\mu\hat{r}(\s_i)$
		\\ \hline
		capped-$\ell_1$ & $ \mu\min(\s_i, \theta ), \; \theta > 0$ \\ \hline
		LSP & $\mu\log \left(\frac{\s_i}{\theta} + 1\right), \; \theta > 0$  \\ \hline
		TNN 
		& $\begin{cases}
		\mu\s_i & i  > \theta \\
		0         & i \le \theta
		\end{cases}$
		\\ \hline
		SCAD & 
		$\begin{cases}
		\mu \s_i                                                  & \sigma_i \le \mu             \\
		\frac{-\s_i^2 + 2 \theta \mu \s_i - \mu^2}{2(\theta - 1)} & \mu < \sigma_i \le \theta\mu \\
		\frac{(\theta + 1) \mu^2}{2}                              & \sigma_i > \theta\mu
		\end{cases}, \theta > 2 $
		\\ \hline
		MCP & 
		$\begin{cases}
		\mu \s_i - \frac{\s_i^2}{2 \theta} & \s_i \le \theta \mu\\
		\frac{\theta \mu^2}{2}             & \s_i >   \theta \mu
		\end{cases}, \theta > 0$
		\\ \hline
	\end{tabular}
	\label{tab:lwregu}
\end{table}

The Iteratively Reweighted Nuclear Norm (IRNN) algorithm \cite{canyi2014}
is a state-of-the-art solver for nonconvex low-rank minimization. It
is based on upper-bounding the nonconvex $r$, and
approximates the matrix rank by 
a weighted version of the nuclear norm
$\| X \|_w = \sum_{i = 1}^m w_i \sigma_i$, 
where $0 \leq w_1 \le \dots \le w_m$, 
Intuitively, $\| X \|_w$ 
imposes a smaller
penalty on the
larger 
(and more informative)
singular values.
Other solvers that are designed only for specific nonconvex low-rank
regularizers include
\cite{sun2013robust}
(for the capped-$\ell_1$),
\cite{hu2013fast}
(for the TNN),
and 
\cite{wang2013nonconvex} (for  the
MCP).
All these (including IRNN) need SVD in each iteration.
It takes $O(m^2 n)$ time,
and thus can be slow.

While proximal gradient algorithms have mostly been used  on convex problems,
recently they are also applied to nonconvex ones \cite{sun2013robust,hu2013fast,canyi2014,lu2015generalized}.
In particular, in the very
recent generalized proximal gradient (GPG) algorithm
\cite{lu2015generalized},
it is shown that
for any nonconvex $r$
satisfying assumption~A3, 
its proximal operator
can be computed by the following generalized 
singular value thresholding
(GSVT) operator.

\begin{prop} \label{pr:proxReduce} \cite{lu2015generalized}
	$\Prox{\mu r}{X} = U\Diag(y^*)V^{\top}$, where $U \Sigma V^{\top}$ is the SVD of $X$, and
	$y^*= [y^*_i]$ with
	\begin{equation} \label{eq:proRed}
	y_i^* \in \arg\min_{y_i \ge 0} \frac{1}{2} \left(y_i - \sigma_i\right)^2 + \mu
	\hat{r}(y_i).
	\end{equation}
\end{prop}

In GPG, problem (\ref{eq:proRed}) is solved by a fixed-point iteration algorithm.
Indeed,
closed-form solutions exist
for the regularizers in Table~\ref{tab:lwregu}
\cite{gong2013general}.
While the obtained proximal operator can then be immediately plugged into a proximal gradient algorithm,
Proposition~\ref{pr:proxReduce}
still involves SVD.

%-----------------------------------------------------------------------------------------
% Proximal Algorithm
%-----------------------------------------------------------------------------------------

\section{Proposed Algorithm}
\label{sec:proalg}

In this section, we show that the GSVT operator $\Prox{\mu r}{\cdot}$ can be computed more efficiently.
It is based on two ideas. First, 
the singular values
in $\Prox{\mu r}{\cdot}$ 
are automatically thresholded. Second,
$\Prox{\mu r}{\cdot}$
can be obtained from the proximal operator on a smaller matrix.

%%%%%%%%%%%%%%%%

\subsection{Automatic Thresholding of Singular Values}

The following Proposition shows that $y_i^*$
in (\ref{eq:proRed}) becomes zero when
$\s_i$ is smaller than a regularizer-specific threshold. Because of the lack of
space, proofs will be reported in a longer version of this paper.

\begin{prop} 
	\label{pr:proxSolution} 
	For any $\hat{r}$ satisfying Assumption A3, there exists a threshold $\gamma > 0$
	such that 
	$y^*_i = 0$
	when $\sigma_i \le \gamma$.
\end{prop} 

By examining the optimality conditions of (\ref{eq:proRed}),
simple closed-form solutions 
can be obtained for 
the nonconvex regularizers in Table~\ref{tab:lwregu}.

\begin{corollary} 
	\label{cor:proxSolution} 
	For the nonconvex regularizers in Table~\ref{tab:lwregu}, 
	their $\gamma$   values are equal to
	\begin{itemize}
		\item capped-$\ell_1$: $\gamma = \min\left(\mu, \theta + \frac{\mu}{2}\right)$;
		\item LSP: $\gamma=\min\left(\frac{\mu}{\theta},\theta\right)$; 
		\item TNN: $\gamma = \max\left(\mu, \sigma_{\theta + 1}\right)$;
		\item SCAD: $\gamma = \mu$;
		\item MCP: $\gamma = \sqrt{\theta} \mu$ if $0 < \theta < 1$, and $\mu$ otherwise.
	\end{itemize}
\end{corollary} 

\begin{algorithm}[ht]
	\caption{Power method to obtain an approximate left subspace of $Z$.}
	\begin{algorithmic}[1]
		\REQUIRE matrix $Z \in \R^{m \times n}$, $R \in \R^{n \times k}$.
		\STATE $Y^1 \leftarrow Z R$;
		\FOR{$t = 1,2,\dots,T_{pm} $}
		\STATE $Q^{t + 1} = \text{QR} (Y^t)$; \quad // QR decomposition
		\STATE $Y^{t + 1} = Z({Z}^{\top} Q^{t + 1})$;
		\ENDFOR
		\STATE return $Q^{T_{pm} + 1}$. 
	\end{algorithmic}
	\label{alg:powermethod}
\end{algorithm}

Proposition~\ref{pr:proxSolution} 
suggests that in each proximal iteration $t$, we only need to compute the leading singular
values/vectors of the matrix iterate $Z^t$.
The power method
(Algorithm~\ref{alg:powermethod})
\cite{halko2011finding} 
is a fast and accurate
algorithm for obtaining an approximation 
of such a subspace. Besides the power method,
algorithms such as PROPACK \cite{larsen1998lanczos} have also been
used \cite{toh2010accelerated}. However, the power method is
more efficient than PROPACK \cite{halko2011finding}.
It
also allows warm-start, which is particularly useful because
of the iterative nature of the proximal gradient algorithm.

%%%%%%%%%%%%%%%%

\subsection{Proximal Operator on a Smaller Matrix}

Assume that $Z^t$ has $\hat{k} \le n$ singular values larger than $\gamma$, and 
its rank-$\hat{k}$ SVD is $U_{\hat{k}} \Sigma_{\hat{k}} V_{\hat{k}}^{\top}$.
The following Proposition shows that 
$\Prox{\mu r}{Z^t}$
can be obtained from the proximal operator on a smaller matrix.

\begin{prop} \label{pr:approGSVT} 
Assume that $Q \in \R^{m \times k}$, where $k \ge \hat{k}$, is orthogonal and $\Span{U_{\hat{k}}} \subseteq \Span{Q}$. Then,
$\Prox{\mu r}{Z^t} = Q \cdot \Prox{\mu r}{Q^{\top} Z^t}$.
\end{prop}

Though SVD is still needed to obtain $\Prox{\mu r}{Q^{\top} Z^t}$,
$Q^{\top} Z^t$ is much smaller than $Z^t$ ($k \times n$ vs $m \times n$).
This smaller SVD takes
$O(n k^2)$ time, 
and the other matrix multiplication steps take $O(m n k)$ time.
Thus, the time complexity for this SVD step is reduced from $O(m^2 n)$ to $O((m + k) n k)$.

%%%%%%%%%%%%%%%%

\subsection{Complete Procedure}
\label{sec:procedure}

The complete procedure 
(Algorithm~\ref{alg:FaNCL})
will be called
FaNCL (\underline{Fa}st \underline{N}on\underline{C}onvex \underline{L}owrank).
The core steps are 
9--16.
We first use the power method to
efficiently obtain an approximate $Q$, whose singular values are then thresholded according to 
Corollary~\ref{cor:proxSolution}.
With $k \ge \hat{k}$, 
the  rank of
$\tilde{X}^p$ will be equal to that of 
$\Prox{\mu r}{Z^t}$.
In each iteration,
we ensure a sufficient decrease of the objective:
\begin{align}
F(X^{t+1}) \le F(X^t)- c_1 \FN{X^{t+1}-X^{t}}^2,
\label{eq:decrease}
\end{align}
where $c_1 = \frac{\tau - \rho}{4}$;
otherwise, the power method is restarted.
Moreover, similar to \cite{hsieh2014nuclear,toh2010accelerated},
steps~6-7 use the column spaces of the previous iterates ($V^{t}$ and $V^{t - 1}$)
to warm-start the power method.
For further speedup, we employ a
continuation strategy
as in \cite{mazumder2010spectral,canyi2014,toh2010accelerated}. Specifically,
$\lambda^t$ is initialized to a large value
and then
decreases
gradually.

\begin{algorithm}[ht]
	\caption{FaNCL (Fast NonConvex Low-rank).}
	\label{alg:FaNCL}
	\begin{algorithmic}[1]
		\STATE choose $\tau > \rho$, $c_1 = \frac{\tau - \rho}{4}$, $\lambda^0 > \lambda$ and $\nu \in (0,1)$;
		\STATE initialize $V_0, V_1 \in \R^{n \times k}$ as random Gaussian matrices, and $X^1 = 0$;
		\FOR{$t = 1,2,\dots T$}
		\STATE $\lambda^t \la (\lambda^{t-1} -\lambda) \nu + \lambda$;
		\STATE $Z^t \la X^t - \frac{1}{\tau} \nabla \f{X^t}$;
		\STATE $V^{t-1} \la V^{t-1} - V^t ({V^t}^{\top}V^{t-1})$, and \\
		remove any zero columns;
		\STATE $R^1 \la \text{QR}([V^{t}, V^{t-1}])$;
		\FOR{$p = 1,2,\dots$}
		\STATE $Q \la \text{PowerMethod}(Z^t, R^p)$;
		\STATE $[ U_A^p, \Sigma_A^p, V_A^p ] \la \text{SVD}(Q^{\top}Z^t)$;
		\STATE $\hat{k} \leftarrow $ number of $\sigma_A$'s are $>\gamma$
		in Corollary~\ref{cor:proxSolution};
		\STATE $\tilde{U}^p \la $ $\hat{k}$ leading columns of $U^p_A$;
		\STATE $\tilde{V}^p \la $ $\hat{k}$ leading columns of $V^p_A$;
		\FOR{$i = 1,2,\dots,\hat{k}$}
		\STATE obtain $y^*_i$ from (\ref{eq:proRed}) with $\mu = 1/\tau$ and $\lambda^t$;
		\ENDFOR
		\STATE $\tilde{X}^p \la (  Q \tilde{U}^p ) \Diag(y_1^*,\dots,y_{\hat{k}}^*) (\tilde{V}^p)^{\top} $;
		\IF{$F(\tilde{X}^p) \le F( X^{t}) - c_1 \|\tilde{X}^p - X^{t}\|_F^2$}
		\STATE $X^{t + 1} \la \tilde{X}^p$, \; $V^{t + 1} \la \tilde{V}^p$;
		\STATE break;
		\ELSE  
		\STATE $R^{p+1} = V_A^p$;
		\ENDIF
		\ENDFOR
		\ENDFOR 
		\RETURN $X^{T + 1}$.
	\end{algorithmic}
\end{algorithm}

Algorithm~\ref{alg:FaNCL} can also be used with the 
nuclear norm. It can be shown that 
the threshold $\gamma$ is equal to $ \lambda / \tau$,
and $y_i^*$ in step~15 has the closed-form solution $\max(\sigma_i - \lambda^t / \tau, 0)$.
However, since our focus is on 
nonconvex regularizers, 
using Algorithm~\ref{alg:FaNCL} 
for nuclear norm minimization will not be further pursued in the sequel.

The power method
%being very useful, 
has also been recently used 
to approximate the SVT in nuclear norm minimization 
\cite{hsieh2014nuclear}.
However, 
\cite{hsieh2014nuclear}
is based on active  subspace selection (which
uses SVT to update the active row and column subspaces
of the current solution), and is thus very different from the proposed algorithm
(which is a proximal gradient algorithm).
In Section~\ref{sec:expts}, it will be shown that the proposed method has better empirical
performance.
Moreover, 
\cite{hsieh2014nuclear}
is only designed for nuclear norm minimization,
and cannot be extended for the nonconvex regularizers considered here.

A breakdown of the time complexity of Algorithm~\ref{alg:FaNCL} is as follows.
For simplicity, assume that $X^t$'s always have rank $k$.  
Step~5 takes $O(mn)$ time;
step~6 and 7 take $O(n k^2)$ time;
step~9 and 10 take $O(mnk T_{pm})$ time;  step~17 takes $O(mnk)$ time;
and step~18 takes  $O(mn)$ time.  Thus, the per-iteration time complexity 
%of Algorithm~\ref{alg:FaNCL} 
is $O(mnk p T_{pm})$.
In the experiment, we set $T_{pm} = 3$ and $p = 1$. Empirically, this setting is enough to guarantee (\ref{eq:decrease}).
In contrast, SVDs in GPG and IRNN take $O(m^2 n)$ time, and are thus much slower as $k \ll m$.

%%%%%%%%%%%%%%%%%%%%%%%%%%%%%%%%%%%%%%%%%%%%%%%%

\subsection{Convergence Analysis}
\label{sec:convergence}

The following Proposition shows that $\{X^t\}$ 
from Algorithm~\ref{alg:FaNCL}
converges to a limit point $X^{\infty} = \lim_{t \rightarrow \infty} X^t$.

\begin{prop} 
	\label{pr:bound}
	$\sum_{t = 1}^{\infty} \FN{X^{t + 1} - X^t}^2 < \infty$.
\end{prop}

%\begin{proof}
%	Summing (\ref{eq:decrease}) from $t = 1$ to $\infty$, we have
%	$F(X^{1}) - F(X^{\infty}) \ge c_1 \sum_{t = 1}^{\infty} \FN{X^{t+1}-X^{t}}^2$.
%	Result then follows from assumption (A2).
%\end{proof}

The following shows that it is also a critical point.\footnote{Since $r$ is nonconvex and its subdifferential for points in its domain may be empty, we 
define $X^*$ as a critical point by extending the definition in
\cite{gong2013general}, namely that
$0\in \nabla f(X^*) + \lambda \partial r_1(X^*) - \lambda \partial r_2(X^*)$, where
$r(X)=r_1(X) - r_2(X)$, and $r_1$ and $r_2$ are convex.}

\begin{theorem} \label{the:convergence}
	$\{X^t\}$ converges to a critical point
	$X^*$ of problem~(\ref{eq:problem})
	in a finite number of iterations.
\end{theorem}

By combining with Proposition~\ref{pr:bound}, 
the following shows that $\FN{X^{t + 1} - X^t}^2$ converges to zero at a rate of $O(1/T)$.

\begin{corollary} 
	\label{the:conv:rate}
	$\min_{t = 1, \dots, T} \FN{X^{t + 1} - X^t}^2 \le \frac{1}{c_1 T} \left[  F(X^1) - F(X^*) \right]$.
\end{corollary}

%-----------------------------------------------------------------------------------------
% Special Formulations
%-----------------------------------------------------------------------------------------

%%%%%%%%%%%%%%%%%%%%%%%%%%%%%%%%%%%%%%%%%%%%%%%%

\subsection{Further Speedup for Matrix Completion}
\label{sec:speedupmatcomp}

In matrix completion, one attempts to recover 
a low-rank matrix $O \in \R^{m \times n}$ 
by 
observing only some of its elements.
Let the 
observed positions be indicated by  
$\Omega \in \{0,1\}^{m \times n}$, such that
$\Omega_{ij}=1$ if $O_{ij}$ is observed, and 0 otherwise.
It can be formulated as an optimization problem
in
(\ref{eq:problem}),  with
$f(X) =\frac{1}{2}\|\SO{X - O}\|_F^2$, where 
%$\SO{\cdot}$ is defined as
$[\SO{A}]_{ij} = 
A_{ij}$ if  $\Omega_{ij} = 1$ and 0      otherwise,
and $r$ is a low-rank regularizer.

It can be easily seen
that step~5
in Algorithm~\ref{alg:FaNCL}
becomes
$Z^t = X^t - \frac{1}{\tau}\SO{X^t - O}$.
By observing that $X^t$ is low-rank and $\frac{1}{\tau}\SO{X^t - O}$ is sparse,
Mazumder \etal
\cite{mazumder2010spectral} showed that
this ``sparse plus low-rank'' structure allows 
matrix multiplications of the form $Z A$ and $Z^{\top}B$ 
%(where $A \in \R^{n \times k}$ and $B \in \R^{m \times k}$)
to be efficiently  computed.
Here, this trick can also be directly used to speed up the computation of $Z({Z}^{\top} Q^{t + 1})$ in 
Algorithm~\ref{alg:powermethod}. 
Since $\|\Omega\|_1$ is very sparse,
this step takes $O(k T_{pm} \|\Omega\|_1)$ time instead of $O( m n k T_{pm})$,
thus is much faster.

The following Proposition shows that
$\FN{ \tilde{X}^p - X^t}^2$ in step~18 
of Algorithm~\ref{alg:FaNCL} can 
also 
be 
easily computed.

\begin{prop} \label{pr:checkobject}
	Let the reduced SVD of $X$ be $U \Sigma V^{\top}$, 
	and $P, Q$ be orthogonal matrices such that $\Span{U} \subseteq \Span{P}$ and $\Span{V} \subseteq \Span{Q}$. 
	Then $\FN{X} = \FN{P^{\top} X Q}$.
\end{prop}
%Proof is in Appendix~\ref{app:checkobject}. 
Let the reduced SVDs of $\tilde{X}^p$ and $X^t$ be $\tU \tS\tV^{\top}$ and $U^t \S^t {V^t}^{\top}$, respectively.
Let $P = \text{QR}([\tU,U^t])$ and $Q = \text{QR}([\tV,V^t])$. 
Using Proposition~\ref{pr:checkobject},
%\begin{eqnarray*}
$\FN{ \tilde{X}^p - X^t} = \FN{P^{\top}(\tilde{X}^p - X^t)Q} 
	= \FN{( P^\top\tU) \tS ( \tV^\top Q ) - (P^\top U^t) \Sigma^t ({V^t}^{\top}Q)}$.
%\end{eqnarray*}
This takes $O(n k^2)$ 
instead of $O(mn)$
time.
The per-iteration time complexity
is reduced from $O(mnk T_{pm})$ to $O( (nk + T_{pm} |\Omega|_1) k )$ and is much faster.
Table~\ref{tab:timecomp}
compares the per-iteration time complexities and convergence rates for the various
low-rank matrix completion solvers 
used in the experiments (Section~\ref{sec:expt1}).

\begin{table}[ht]
	\centering
	\renewcommand{\arraystretch}{1.2}
	\caption{Comparison of the per-iteration time complexities and convergence rates of various matrix completion solvers. 
		Here, $\nu \in (0,1)$ is a constant.}
	\vspace{-5px}
	\begin{tabular}{c | c | c | c}
		\hline
		 regularizer  & method                                          & complexity                 & rate         \\ \hline
		  (convex)    & APG \cite{ji2009accelerated,toh2010accelerated} & $O(m n k)$                 & $O(1/T^2)$   \\ \cline{2-4}
		   nuclear    & Soft-Impute \cite{mazumder2010spectral}         & $O(k\|\Omega\|_1)$         & $O(1/T)$     \\ \cline{2-4}
		    norm      & active ALT \cite{hsieh2014nuclear}              & $O(k T_{in} \|\Omega\|_1)$ & $O(\nu^{T})$ \\ \hline
		 fixed-rank   & LMaFit \cite{wen2012solving}                    & $O(k \|\Omega\|_1)$        & ---          \\ \cline{2-4}
		factorization & R1MP \cite{wang2014rank}                        & $O(\|\Omega\|_1)$          & $O(\nu^T)$   \\ \hline
		  nonconvex   & IRNN \cite{canyi2014}                           & $O(m^2 n)$                 & ---          \\ \cline{2-4}
		              & GPG \cite{lu2015generalized}                    & $O(m^2 n)$                 & ---          \\ \cline{2-4}
		              & FaNCL                                           & $O(k \|\Omega\|_1)$        & $O(1/T)$     \\ \hline
	\end{tabular}
	\label{tab:timecomp}
\end{table}

%%%%%%%%%%%%%%%%%%%%%%%%%%%%%%%%%%%%%%%%%%%%%%%%

\subsection{Handling Multiple Matrix Variables}

%Besides solving (\ref{eq:problem}), 
The proposed algorithm can be extended for 
optimization problems 
involving
$N$ matrices $X_1, \dots, X_N$:
\begin{align} 
\!\!\!\!\!\!\!
\min
%\min_{X_i}  
F(X_1, \dots, X_N) \! \equiv\! 
f(X_1, \dots, X_N) \! + \! \sum_{i = 1}^N \lambda_i r_i(X_i).
\label{eq:problem:gen}
\end{align}
Assumptions A1-A3 are analogously extended. In
particular, 
A1 now assumes that
$\FN{\nabla f_i(X) - \nabla f_i(Y)} \le \rho_i \FN{X - Y}$ for some $\rho_i$,
where $f_i(X)$ is the function obtained by keeping all the $X_j$'s (where $i\ne j$) in $f$
fixed.

Many machine learning problems can be cast into this form. One example that will be considered
in Section~\ref{sec:expts} is 
robust principal component analysis (RPCA) \cite{candes2011robust}.
Given a noisy data matrix $O$, RPCA assumes that $O$ can be approximated by the sum of a
low-rank matrix $X$
plus sparse data noise $Y$. 
Mathematically, we have 
%the following optimization problem:
\begin{equation} 
\label{eq:rpca}
\min_{X,Y} F(X,Y) \equiv f(X,Y) + \lambda r(X) + \beta \|Y\|_1,
\end{equation} 
where 
$f(X, Y) = \frac{1}{2}\|X + Y - O\|_F^2$,
$r$ is a low-rank regularizer on $X$, and $\|Y\|_1$ encourages $Y$ to be sparse.
Since both $r$ and the $\ell_1$ regularizer $\|\cdot\|_1$ are nonsmooth, 
(\ref{eq:rpca}) does not fit into formulation (\ref{eq:problem}). 
Besides RPCA, 
problems such as subspace clustering \cite{liu2013robust}, 
multilabel learning \cite{zhu2010image} and multitask learning
\cite{gong2012robust}
can also be cast as (\ref{eq:problem:gen}).

For simplicity, we focus on the 
case with
two parameter blocks.  Extension to multiple blocks
is straightforward.
To solve the two-block problem
in (\ref{eq:rpca}),
we perform alternating proximal steps on $X$ and $Y$ at each iteration $t$:
\begin{eqnarray*}
	X^{t + 1}
	= & \!\!\!\! \arg\min_X \frac{1}{2} \FN{X- Z_X^t}^2 + \frac{\lambda}{\tau}r(X)
	= & \!\!\!\! \Prox{\frac{\lambda}{\tau}r}{Z^t_X},\\
	Y^{t + 1} = & \!\!\!\! \arg\min_Y \frac{1}{2} \FN{Y- Z_Y^t}^2 + \frac{\beta}{\tau}\|Y\|_1
	= & \!\!\!\! \Prox{\frac{\beta}{\tau}\|\cdot\|_1}{Z^t_Y},
\end{eqnarray*}
where 
$Z_X^t = X^t - \frac{1}{\tau} \nabla \f{X^t,Y^t}$, 
and 
$Z_Y^t = Y^t - \frac{1}{\tau} \nabla \f{X^{t+1},Y^t}$. 
$Y^{t + 1}$ can be easily obtained as 
$Y^{t + 1}_{ij} =
\text{sign} \left([Z_Y^t]_{ij} \right) \left( \left|[Z_Y^t]_{ij} \right| - \frac{\beta}{\tau}
\right)_+$,
where $\text{sign}(x)$ denotes the sign of $x$.
Similar to (\ref{eq:decrease}), we ensure a sufficient decrease of the objective 
in each iteration: 
\begin{eqnarray*}
	F_{Y^t}(X^{t+1}) & \le  & F_{Y^t}(X^t) - c_1 \FN{X^{t+1} - X^{t}}^2,\\
	F_{X^{t+1}}(Y^{t+1})  & \le & F_{X^{t+1}}(Y^t) - c_1 \FN{Y^{t+1} - Y^{t}}^2,
\end{eqnarray*}
where 
$F_Y(X) = f(X, Y)+ \lambda r(X)$, 
and 
$F_X(Y) = f(X, Y) + \beta \|Y\|_1$.
The resultant algorithm is similar to Algorithm~\ref{alg:FaNCL}.

When $F$ is convex,  convergence of this 
alternating minimization scheme has been well studied \cite{tseng2001convergence}.
However, here $F$ is nonconvex. 
We extend
the convergence results in Section~\ref{sec:convergence} 
to the following.

\begin{theorem} 
	\label{the:conv:altprox}
	With $N$ parameter blocks and $\{(X_1^t, \dots, X_N^t)\}$ generated by the algorithm, we have
	\begin{enumerate}
		\item $\sum_{t = 1}^{\infty} \sum_{i = 1}^N \FN{X_i^{t+1}-X_i^t}^2 < \infty$;
		\item 
		$\{(X_1^t, \dots, X_N^t)\}$ 
		converges to a critical point $(X_1^*, \dots, X_N^*)$ of (\ref{eq:problem:gen}) in a finite
		number of iterations;
		\item $\min_{t = 1, \dots, T} \sum_{i = 1}^N \FN{X_i^{t + 1} - X_i^t}^2 \le 
		\frac{1}{c_1 T} [ F(X_1^1, \dots, X_N^1) - F(X_1^*, \dots, X_N^*) ]$.
	\end{enumerate}
\end{theorem}

%%%%%%%%%%%%%%%%%%%%%%%%%%%%%%%%%%%%%%%%%%%%%%%%

\section{Experiments}
\label{sec:expts}

%In this section, we perform extensive experiments on matrix completion (Section~\ref{sec:expt1}) and RPCA (Section~\ref{sec:expt2}).

%%%%%%%%%%%%%%%%%%%%%%%%%%%%%%%%%%%%%%%%%%%%%%%%

\subsection{Matrix Completion}
\label{sec:expt1}

\begin{table*}[t]
	\centering
	\caption{Matrix completion performance on the synthetic data.  Here, NMSE
		is scaled by $\times 10^{-2}$, and CPU time is in seconds.}
	\scriptsize
	\vspace{-5px}
	\begin{tabular}{cc|ccc|ccc|ccc|ccc}
		\hline
		         &              &         \multicolumn{3}{c|}{$m=500$}         &        \multicolumn{3}{c|}{$m=1000$}         &        \multicolumn{3}{c|}{$m=1500$}         &         \multicolumn{3}{c}{$m=2000$}         \\
		         &              &  \multicolumn{3}{c|}{(observed: $12.43\%$)}  &  \multicolumn{3}{c|}{(observed: $6.91\%$)}   &  \multicolumn{3}{c|}{(observed: $4.88\%$)}   &   \multicolumn{3}{c}{(observed: $3.80\%$)}   \\
		         &              &           NMSE            & rank &   time    &           NMSE            & rank &   time    &           NMSE            & rank &   time    &           NMSE            & rank &   time    \\ \hline
		nuclear  &     APG      &    $3.95 \!\pm\! 0.16$    &  49  &    4.8    &    $3.90 \!\pm\! 0.05$    &  59  &   59.5    &   $3.74 \!\pm \! 0.02$    &  71  &   469.3   &    $3.69 \!\pm\! 0.04$    &  85  &  1383.3   \\ \cline{2-14}
		  norm   & Soft-Impute  &    $3.95 \!\pm\! 0.16$    &  49  &   64.9    &    $3.90 \!\pm\! 0.05$    &  59  &   176.0   &    $3.74 \!\pm\! 0.02$    &  71  &   464.4   &    $3.69 \!\pm\! 0.04$    &  85  &  1090.2   \\ \cline{2-14}
		         & active   ALT &    $3.95 \!\pm\! 0.16$    &  49  &   17.1    &    $3.90 \!\pm\! 0.05$    &  59  &   81.9    &    $3.74 \!\pm\! 0.02$    &  71  &   343.8   &    $3.69 \!\pm\! 0.04$    &  85  &   860.1   \\ \hline
		 fixed   &    LMaFit    &    $2.63 \!\pm\! 0.10$    &  5   &    0.6    &    $2.85 \!\pm\! 0.10$    &  5   &    1.7    &    $2.54 \!\pm\! 0.09$    &  5   &    4.5    &    $2.40 \!\pm\! 0.09$    &  5   &    7.1    \\ \cline{2-14}
		  rank   &     R1MP     &   $22.72 \!\pm\! 0.63$    &  39  &    0.3    &  $20.89 \!\pm\!
		0.66$   &  54  &    0.8    &   $20.04 \!\pm\! 0.66$    &  62  &    1.4    &  $19.53 \!\pm\!
		0.61$   &  63  &    3.4    \\ \hline\hline
		 capped  &     IRNN     & ${\bf 1.98} \!\pm\! 0.07$ &  5   &    8.5    & ${\bf 1.89} \!\pm\! 0.04$ &  5   &   75.5    & ${\bf 1.81} \!\pm\! 0.02$ &  5   &   510.8   & ${\bf 1.80} \!\pm\! 0.02$ &  5   &  1112.3   \\ \cline{2-14}
		$\ell_1$ &     GPG      & ${\bf 1.98} \!\pm\! 0.07$ &  5   &    8.5    & ${\bf 1.89} \!\pm\! 0.04$ &  5   &   72.4    & ${\bf 1.81} \!\pm\! 0.02$ &  5   &   497.0   & ${\bf 1.80} \!\pm\! 0.02$ &  5   &  1105.8   \\ \cline{2-14}
		         &    FaNCL     & ${\bf 1.98} \!\pm\! 0.07$ &  5   & {\bf 0.3} & ${\bf 1.89}\!\pm\! 0.04$  &  5   & {\bf 0.9} & ${\bf 1.81} \!\pm\! 0.02$ &  5   & {\bf 2.6} & ${\bf 1.80} \!\pm\! 0.02$ &  5   & {\bf 4.1} \\ \hline
		  LSP    &     IRNN     & ${\bf 1.98} \!\pm\! 0.07$ &  5   &   21.8    & ${\bf 1.89} \!\pm\! 0.04$ &  5   &   223.9   & ${\bf 1.81} \!\pm\! 0.02$ &  5   &   720.9   & ${\bf 1.80} \!\pm\! 0.02$ &  5   &  2635.0   \\ \cline{2-14}
		         &     GPG      & ${\bf 1.98} \!\pm\! 0.07$ &  5   &   21.2    & ${\bf 1.89} \!\pm\! 0.04$ &  5   &   235.3   & ${\bf 1.81} \!\pm\! 0.02$ &  5   &   687.4   & ${\bf 1.80} \!\pm\! 0.02$ &  5   &  2612.0   \\ \cline{2-14}
		         &    FaNCL     & ${\bf 1.98} \!\pm\! 0.07$ &  5   & {\bf 0.5} & ${\bf 1.89} \!\pm\! 0.04$ &  5   & {\bf 2.2} & ${\bf 1.81} \!\pm\! 0.02$ &  5   & {\bf 3.3} & ${\bf 1.80} \!\pm\! 0.02$ &  5   & {\bf 7.6} \\ \hline
		  TNN    &     IRNN     & ${\bf 1.98} \!\pm\! 0.07$ &  5   &    8.5    & ${\bf 1.89} \!\pm\! 0.04$ &  5   &   72.6    & ${\bf 1.81} \!\pm\! 0.02$ &  5   &   650.7   & ${\bf 1.80} \!\pm\! 0.02$ &  5   &  1104.1   \\ \cline{2-14}
		         &     GPG      & ${\bf 1.98} \!\pm\! 0.07$ &  5   &    8.3    & ${\bf 1.89} \!\pm\! 0.04$ &  5   &   71.7    & ${\bf 1.81} \!\pm\! 0.02$ &  5   &   655.3   & ${\bf 1.80} \!\pm\! 0.02$ &  5   &  1098.2   \\ \cline{2-14}
		         &    FaNCL     & ${\bf 1.98} \!\pm\! 0.07$ &  5   & {\bf 0.3} & ${\bf 1.89} \!\pm\! 0.04$ &  5   & {\bf 0.8} & ${\bf 1.81} \!\pm\! 0.02$ &  5   & {\bf 2.7} & ${\bf 1.80} \!\pm\! 0.02$ &  5   & {\bf 4.2} \\ \hline
	\end{tabular}
	\label{tab:sythmatcomp}
\end{table*}

We compare a number of low-rank matrix completion solvers, including
models based on 
(i) the commonly used (convex) nuclear norm regularizer; 
(ii) fixed-rank factorization models \cite{wen2012solving,wang2014rank},
which decompose the  observed matrix $O$ into a product of 
rank-$k$ matrices $U$ and $V$. Its optimization problem can be written as:
$\min_{U,V}\frac{1}{2}\FN{\SO{U V - O}}^2+ \frac{\lambda}{2}(\|U\|_F^2 + \|V\|_F^2)$;
and 
(iii) nonconvex regularizers, including the capped-$\ell_1$ (with $\theta$ in
Table~\ref{tab:lwregu} set to $2 \lambda$), LSP
(with $\theta = \sqrt{\lambda}$), and TNN
(with $\theta = 3$).

The nuclear norm minimization algorithms to be compared include:
\begin{enumerate}
	\item Accelerated proximal gradient
	(APG)\footnote{\url{http://perception.csl.illinois.edu/matrix-rank/Files/apg_partial.zip}}
	algorithm 
	\cite{ji2009accelerated,toh2010accelerated}, with
	the partial SVD by PROPACK \cite{larsen1998lanczos};
	
	\item Soft-Impute\footnote{\url{http://cran.r-project.org/web/packages/softImpute/index.html}}
	\cite{mazumder2010spectral}, which iteratively replaces the missing elements with those obtained
	from SVT. The ``sparse plus low-rank'' structure of  the
	matrix iterate is utilized to speed up computation (Section~\ref{sec:speedupmatcomp});
	
	\item Active alternating minimization\footnote{\url{http://www.cs.utexas.edu/~cjhsieh/nuclear_active_1.1.zip}}
	(denoted ``active ALT'') \cite{hsieh2014nuclear}, which adds/removes rank-one subspaces from the active set in each iteration.  The nuclear norm optimization problem is 
	then reduced to a smaller problem defined only on this active set.
\end{enumerate}

We do not compare
with the Frank-Wolfe algorithm \cite{zhang2012accelerated} and stochastic gradient descent \cite{avron2012efficient},
as they have been shown to be less efficient 
\cite{hsieh2014nuclear}.
For the fixed-rank factorization models 
(where the rank is tuned by the validation set),
we compare with the two state-of-the-art algorithms:
\begin{enumerate}
	\item Low-rank matrix fitting (LMaFit)
	algorithm\footnote{\url{http://www.caam.rice.edu/~optimization/L1/LMaFit/download.html}}
	\cite{wen2012solving}; and
	
	\item Rank-one matrix pursuit (R1MP) \cite{wang2014rank},
	which pursues a rank-one basis in each iteration.
\end{enumerate}
We do not compare with the concave-convex procedure \cite{zhang2010analysis,hu2013fast}, since
it has been shown to be inferior to IRNN \cite{gong2013general}.

For models with nonconvex low-rank regularizers, we compare the following solvers:
\begin{enumerate}		
	\item Iterative reweighted nuclear norm (IRNN)\footnote{\url{https://sites.google.com/site/canyilu/file/2014-CVPR-IRNN.zip?attredirects=0&d=1}}
	\cite{canyi2014};
	\item Generalized proximal gradient (GPG) algorithm \cite{lu2015generalized}, with the
	underlying problem (\ref{eq:proRed}) solved more efficiently using the closed-form solutions
	in \cite{gong2013general};
	\item The proposed FaNCL algorithm ($T_{pm} = 3$, $p = 1$).
\end{enumerate}

All algorithms are implemented in Matlab.
The 
same 
stopping criterion is
used,
namely that the algorithm stops when the difference
in objective values between consecutive iterations is smaller than a given
threshold.
Experiments are run on a PC with i7 4GHz CPU and 24GB memory.

%%%%%%%%%%%%%%%

\subsubsection{Synthetic Data}
\label{sec:matcomp:syn}

The observed $m\times m$ matrix is generated as $O = 
U V + G$,  where 
the elements 
of $U \in \R^{m \times k}, V \in \R^{k \times m}$ 
(with $k = 5$) are 
sampled i.i.d. from the normal distribution $\mathcal{N}(0,
1)$, and
elements of $G$ sampled from $\mathcal{N}(0, 0.1)$.
A total of $\|\Omega\|_1 = 2 m k \log(m)$ random elements
in $O$ are observed.  
Half of them are used for training, and the rest as validation
set for parameter tuning.
Testing is performed on the non-observed (missing) elements.

For performance evaluation,  
we use (i) the normalized mean squared error
$\text{NMSE} = \sqrt{\sum_{(i,j) \not\in \Omega} (X_{ij} - [UV]_{ij})^2}/\sqrt{\sum_{(i,j) \not\in \Omega} [UV]_{ij}^2}$,
where
$X$ is
the recovered matrix;
(ii) rank of $X$;
and 
(iii) training CPU time.
We vary $m$ in the range $\{500, 1000, 1500, 2000\}$.
Each experiment is repeated five times.

\begin{figure}[ht]
	\centering
	\includegraphics[width = 0.3\textwidth]{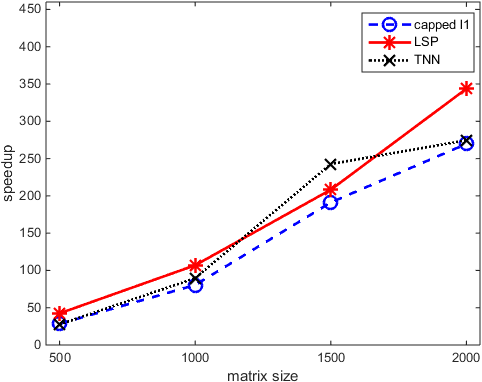}
	\vspace{-.1in}
	\caption{Speedup of FaNCL over GPG at different matrix sizes.}
	\label{fig:speedupmatcomp}
\end{figure}

Results
are shown in Table~\ref{tab:sythmatcomp}.  As can be seen, the nonconvex regularizers (capped-$\ell_1$, LSP and TNN) lead to much lower NMSE's than 
the convex nuclear norm regularizer and fixed-rank factorization.
Moreover, 
as is also observed in \cite{avron2012efficient},
the nuclear norm needs to use a much higher rank than the
nonconvex ones.  
In terms of speed, 
FaNCL is the fastest among the nonconvex low-rank solvers.
Figure~\ref{fig:speedupmatcomp} shows
its speedup over GPG (which in turn is faster than IRNN).
As can be seen, the larger the matrix, the higher is the speedup.
%which can be over 300 times.

Recall that the efficiency of the proposed algorithm comes from (i) automatic singular
value thresholding; (ii) computing the proximal operator on a smaller matrix; and (iii)
exploiting the ``sparse plus low-rank'' structure in matrix completion.
Their individual contributions 
are examined in Table~\ref{tab:fancl}.
The baseline is GPG, which uses none of these; while the proposed FaNCL uses all.
As all the variants produce the same solution,
%(but with different speeds),
the obtained NMSE and rank values
are not shown.
As can be seen, tricks (i), (ii) and (iii) lead to average speedups of about 6, 4, and 3, respectively;
and are particularly useful on the large data sets.

\begin{table}[ht]
	\centering
	\footnotesize
	\caption{Effects of the three tricks
		on CPU time (in seconds) using the synthetic data.
		(i) automatic singular
		value thresholding; (ii) computing the proximal operator on a smaller matrix; and (iii)
		``sparse plus low-rank'' structure.}
	\vspace{-5px}
	\begin{tabular}{c|c|c|c|c|c}
		\hline
		         &             solver              &   $500$   &  $1000$   &  $1500$   &  $2000$   \\ \hline
		 capped  &         baseline (GPG)          &    8.5    &   72.4    &   497.0   &  1105.8   \\ \cline{2-6}
		$\ell_1$ &         $\checkmark$ i          &    5.4    &   37.6    &   114.8   &   203.7   \\ \cline{2-6}
		         &       $\checkmark$ i, ii        &    0.6    &    3.2    &   11.4    &   25.6    \\ \cline{2-6}
		         & $\checkmark$ i, ii, iii (FaNCL) & {\bf 0.3} & {\bf 0.9} & {\bf 2.6} & {\bf 6.8} \\ \hline
		  LSP    &         baseline (GPG)          &   21.2    &   235.3   &   687.4   &  2612.0   \\ \cline{2-6}
		         &         $\checkmark$ i          &    4.9    &   44.0    &   70.0    &   154.9   \\ \cline{2-6}
		         &       $\checkmark$ i, ii        &    1.0    &    9.7    &   14.8    &   31.1    \\ \cline{2-6}
		         & $\checkmark$ i, ii, iii (FaNCL) & {\bf 0.5} & {\bf 2.2} & {\bf 3.3} & {\bf 8.2} \\ \hline
		  TNN    &         baseline (GPG)          &    8.3    &   71.7    &   655.3   &  1098.2   \\ \cline{2-6}
		         &         $\checkmark$ i          &    5.4    &   32.5    &   122.3   &   194.1   \\ \cline{2-6}
		         &       $\checkmark$ i, ii        &    0.6    &    2.8    &   10.3    &   15.8    \\ \cline{2-6}
		         & $\checkmark$ i, ii, iii (FaNCL) & {\bf 0.3} & {\bf 0.8} & {\bf 2.7} & {\bf 3.3} \\ \hline
	\end{tabular}
	\label{tab:fancl}
\end{table}

%%%%%%%%%%%%%%%%%%%%%%%%%%%%%%

\subsubsection{MovieLens}

\begin{figure*}[t]
	\centering
	\subfigure[capped-$\ell_1$. \label{fig:cap}]
	{\includegraphics[width = 0.25\textwidth]{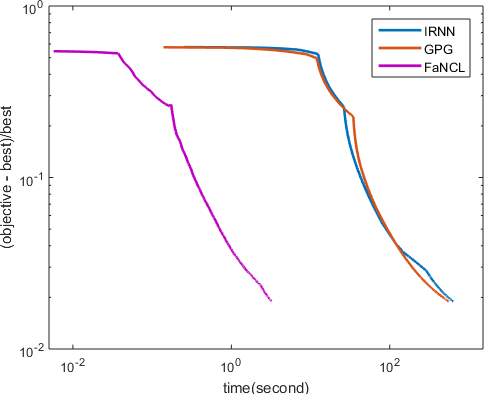}}
	\quad\quad
	\subfigure[LSP. \label{fig:lsp}]
	{\includegraphics[width = 0.25\textwidth]{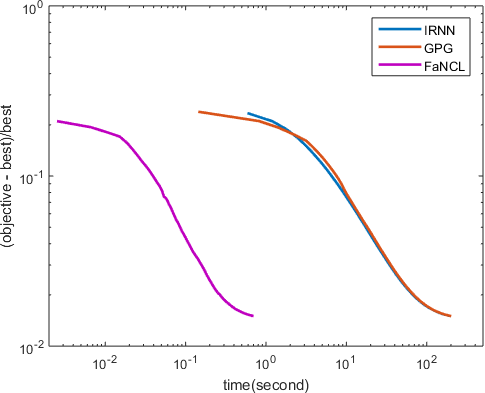}}
	\quad\quad
	\subfigure[TNN. \label{fig:tnn}]
	{\includegraphics[width = 0.25\textwidth]{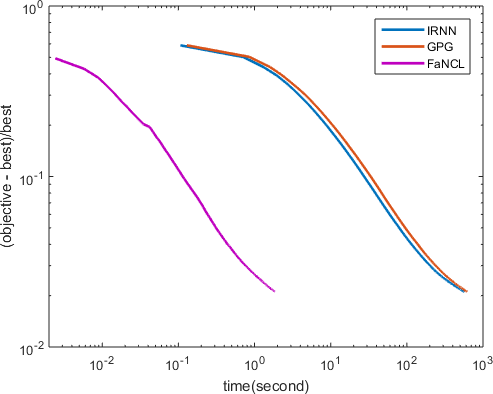}}
	\vspace{-5px}
	\caption{Objective value vs CPU time for the various nonconvex low-rank regularizers on the MovieLens-100K data set.}
	\label{fig:movielen:obj}
\end{figure*}

Experiment is performed on the popular
MovieLens\footnote{\url{http://grouplens.org/datasets/movielens/}} data set
(Table~\ref{tab:recSys}),
which contain ratings of different users on movies.
We follow the setup in \cite{wang2014rank},  
and use $50\%$ of the observed ratings for training, $25\%$ for validation and the rest for testing.
For performance evaluation, we use the root mean squared error on the test set $\Omega$:
$\text{RMSE} = \sqrt{\|\SO{X - O}\|_F^2 / \|\Omega\|_1}$,
where $X$ is the recovered matrix.
The experiment is repeated five times.

\begin{table}[ht]
	\centering
	\caption{Recommendation data sets used in the experiments.}
	\vspace{-5px}
	\begin{tabular}{cc | c | c | c }
		\hline
		          &                  & \#users & \#movies & \#ratings   \\ \hline
		MovieLens &       100K       & 943     & 1,682    & 100,000     \\ \cline{2-5}
		          &        1M        & 6,040   & 3,449    & 999,714     \\ \cline{2-5}
		          &       10M        & 69,878  & 10,677   & 10,000,054  \\ \hline
		\multicolumn{2}{c|}{netflix} & 480,189 & 17,770   & 100,480,507 \\ \hline
		 \multicolumn{2}{c|}{yahoo}  & 249,012 & 296,111  & 62,551,438  \\ \hline
	\end{tabular}
	\label{tab:recSys}
\end{table}

\begin{figure}[H]
	\centering
	\includegraphics[width = 0.3 \textwidth]{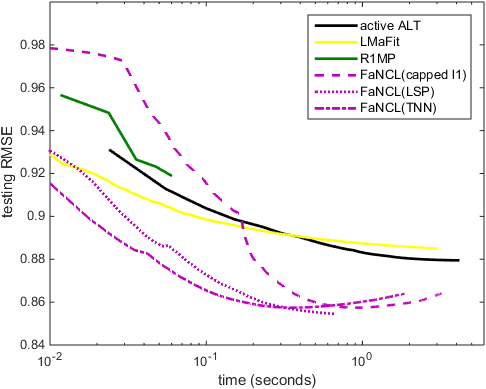}
	\vspace{-.1in}
	\caption{RMSE vs CPU time on the MovieLens-100K data set. }
	\label{fig:movielen:rmse}
	\vspace{-.2in}
\end{figure}

Results are shown in Table~\ref{tab:movielen}. 
Again, nonconvex regularizers lead to the lowest RMSE's. 
Moreover, FaNCL is also the fastest among the nonconvex low-rank solvers,
even faster than the state-of-the-art.
In particular, it is the only solver (among those compared) that can be run on the MovieLens-10M data.
Table~\ref{tab:fancl2} examines the usefulness of the three tricks.
The behavior is similar to
that as observed in 
Table~\ref{tab:fancl}.
Figures~\ref{fig:movielen:obj} and
\ref{fig:movielen:rmse} compare
the objective
and RMSE vs CPU time  for the various methods on the MovieLens-100K data set.
As can be seen,
FaNCL decreases the objective and RMSE much faster than the others.

\begin{table*}[t]
	\centering
	\caption{Matrix completion results on the MovieLens data sets (time is in seconds).}
	\vspace{-5px}
	\label{tab:movielen}
	\begin{tabular}{cc|ccc|ccc|ccc}
		\hline
		                &             &      \multicolumn{3}{c|}{MovieLens-100K}       &        \multicolumn{3}{c|}{MovieLens-1M}        &        \multicolumn{3}{c}{MovieLens-10M}         \\
		                &             &            RMSE             & rank &   time    &            RMSE             & rank &    time    &            RMSE             & rank &    time     \\ \hline
		 nuclear norm   &     APG     &    $0.879 \!\pm\! 0.001$    &  36  &   18.9    &    $0.818 \!\pm\! 0.001$    &  67  &   735.8    &             ---             & ---  & $> 10^{5}$  \\ \cline{2-11}
		                & Soft-Impute &    $0.879 \!\pm\! 0.001$    &  36  &   13.8    &    $0.818 \!\pm\! 0.001$    &  67  &   311.8    &             ---             & ---  & $> 10^{5}$  \\ \cline{2-11}
		                & active ALT  &    $0.879 \!\pm\! 0.001$    &  36  &    4.1    &    $0.818 \!\pm\! 0.001$    &  67  &   133.4    &    $0.813 \!\pm\! 0.001$    & 119  &   3675.2    \\ \hline
		  fixed rank    &   LMaFit    &    $0.884 \!\pm\! 0.001$    &  2   &    3.0    &    $0.818 \!\pm\! 0.001$    &  6   &    39.2    &    $0.795 \!\pm\! 0.001$    &  9   &    650.1    \\ \cline{2-11}
		                &    R1MP     &    $0.924 \!\pm\! 0.003$    &  5   &    0.1    &    $0.862 \!\pm\! 0.004$    &  19  &    2.9     &    $0.850 \!\pm\! 0.008$    &  29  &    37.3     \\ \hline\hline
		capped-$\ell_1$ &    IRNN     &    $0.863 \!\pm\! 0.003$    &  3   &   558.9   &             ---             & ---  & $> 10^{4}$ &             ---             & ---  & $> 10^{5}$  \\ \cline{2-11}
		                &     GPG     &    $0.863 \!\pm\! 0.003$    &  3   &   523.6   &             ---             & ---  & $> 10^{4}$ &             ---             & ---  & $> 10^{5}$  \\ \cline{2-11}
		                &    FaNCL    &    $0.863 \!\pm\! 0.003$    &  3   & {\bf 3.2} &    $0.797 \!\pm\! 0.001$    &  5   & {\bf 29.4} &    $0.783 \!\pm\! 0.002$    &  8   & {\bf 634.6} \\ \hline
		      LSP       &    IRNN     & ${\bf 0.855} \!\pm\! 0.002$ &  2   &   195.9   &             ---             & ---  & $> 10^{4}$ &             ---             & ---  & $> 10^{5}$  \\ \cline{2-11}
		                &     GPG     & ${\bf 0.855} \!\pm\! 0.002$ &  2   &   192.8   &             ---             & ---  & $> 10^{4}$ &             ---             & ---  & $> 10^{5}$  \\ \cline{2-11}
		                &    FaNCL    & ${\bf 0.855} \!\pm\! 0.002$ &  2   & {\bf 0.7} & ${\bf 0.786} \!\pm\! 0.001$ &  5   & {\bf 25.6} & ${\bf 0.777} \!\pm\! 0.001$ &  9   & {\bf 616.3} \\ \hline
		      TNN       &    IRNN     &    $0.862 \!\pm\! 0.003$    &  3   &   621.9   &             ---             & ---  & $> 10^{4}$ &             ---             & ---  & $> 10^{5}$  \\ \cline{2-11}
		                &     GPG     &    $0.862 \!\pm\! 0.003$    &  3   &   572.7   &             ---             & ---  & $> 10^{4}$ &             ---             & ---  & $> 10^{5}$  \\ \cline{2-11}
		                &    FaNCL    &    $0.862 \!\pm\! 0.003$    &  3   & {\bf 1.9} &    $0.797 \!\pm\! 0.004$    &  5   & {\bf 25.8} &    $0.783 \!\pm\! 0.002$    &  8   & {\bf 710.7} \\ \hline
	\end{tabular}
\end{table*}

\begin{table}[ht]
	\centering
	\caption{Effects of the three tricks on CPU time (in seconds) on the MovieLens data.}
	\vspace{-5px}
	\begin{tabular}{c|c|c|c|c}
		\hline
		         &             solver              &   100K    &     1M     &     10M     \\ \hline
		 capped  &         baseline (GPG)          &   523.6   & $> 10^{4}$ & $> 10^{5}$  \\ \cline{2-5}
		$\ell_1$ &         $\checkmark$ i          &   212.2   &   1920.5   & $> 10^{5}$  \\ \cline{2-5}
		         &       $\checkmark$ i, ii        &   29.2    &   288.8    & $> 10^{5}$  \\ \cline{2-5}
		         & $\checkmark$ i, ii, iii (FaNCL) & {\bf 3.2} & {\bf 29.4} & {\bf 634.6} \\ \hline
		  LSP    &         baseline (GPG)          &   192.8   & $> 10^{4}$ & $> 10^{5}$  \\ \cline{2-5}
		         &         $\checkmark$ i          &   35.8    &   2353.8   & $> 10^{5}$  \\ \cline{2-5}
		         &       $\checkmark$ i, ii        &    5.6    &   212.4    & $> 10^{5}$  \\ \cline{2-5}
		         & $\checkmark$ i, ii, iii (FaNCL) & {\bf 0.7} & {\bf 25.6} & {\bf 616.3} \\ \hline
		  TNN    &         baseline (GPG)          &   572.7   & $> 10^{4}$ & $> 10^{5}$  \\ \cline{2-5}
		         &         $\checkmark$ i          &   116.9   &   1944.8   & $> 10^{5}$  \\ \cline{2-5}
		         &       $\checkmark$ i, ii        &   15.4    &   256.1    & $> 10^{5}$  \\ \cline{2-5}
		         & $\checkmark$ i, ii, iii (FaNCL) & {\bf 1.9} & {\bf 25.8} & {\bf 710.7} \\ \hline
	\end{tabular}
	\label{tab:fancl2}
\end{table}

\begin{figure}[H]
	\centering
	\subfigure[netflix.]
	{\includegraphics[width=0.3 \textwidth]{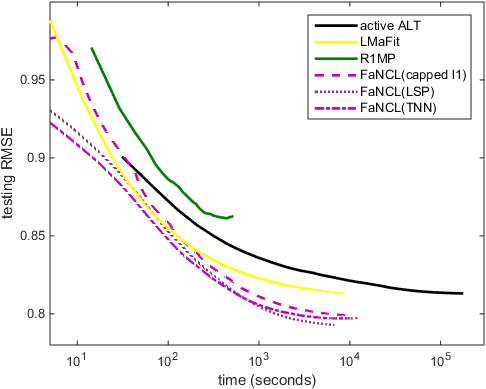}}
	\subfigure[yahoo.]
	{\includegraphics[width=0.3 \textwidth]{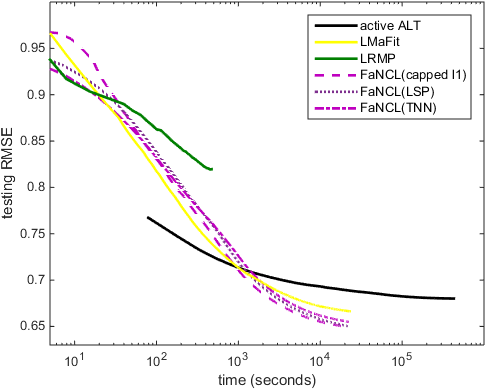}}
	\vspace{-.1in}
	\caption{RMSE vs CPU time on the netflix and yahoo data sets.}
	\label{fig:large:rmse}
\end{figure}

%%%%%%%%%%%%%%%%%%%%%%

\subsubsection{Netflix and Yahoo}

Next, we perform experiments on two very large recommendation
data sets,
Netflix\footnote{\url{http://archive.ics.uci.edu/ml/datasets/Netflix+Prize}} 
and Yahoo\footnote{\url{http://webscope.sandbox.yahoo.com/catalog.php?datatype=c}}
(Table~\ref{tab:recSys}).
We randomly use $50\%$ 
of the observed ratings 
for training,
$25\%$ for validation and the rest for testing.
Each experiment is repeated five times.

Results are shown in Table~\ref{tab:yahoo}.
APG, Soft-Impute, GPG and IRNN cannot be run as the
data set
is large.
Figure~\ref{fig:large:rmse} shows
the objective and RMSE vs time for the compared methods.\footnote{On these two data sets, R1MP easily overfits as the rank increases. Hence, the validation set
	selects a 
	rank 
	which is small 
	(relative to that obtained by the nuclear norm) and R1MP stops
	earlier.
	However, as can be seen, its RMSE is much worse.
}
Again, the nonconvex regularizers converge faster, yield lower RMSE's and 
solutions
of much lower ranks.
Moreover, FaNCL is fast. 

\begin{table*}[htbp]
	\centering
	\caption{Results on the	netflix and yahoo data sets (CPU time is in hours).}
	\vspace{-10px}
	\begin{tabular}{cc|ccc|ccc}
		\hline
		                 &            &     \multicolumn{3}{c|}{netflix}      &       \multicolumn{3}{c}{yahoo}        \\
		                 &            &          RMSE           & rank & time &          RMSE           & rank & time  \\ \hline
		nuclear     norm & active ALT &    $0.814 \pm 0.001$    & 399  & 47.6 &    $0.680 \pm 0.001$    & 221  & 118.9 \\ \hline
		fixed      rank  &   LMaFit   &    $0.813 \pm 0.003$    &  16  & 2.4  &    $0.667 \pm 0.002$    &  10  &  6.6  \\ \cline{2-8}
		                 &    R1MP    &    $0.861 \pm 0.006$    &  31  & 0.2  &   $0.810 \pm
		0.005$   &  92  &  0.3  \\ \hline\hline
		capped-$\ell_1$  &   FaNCL    &    $0.799 \pm 0.001$    &  15  & 2.5  & ${\bf 0.650} \pm 0.001$ &  8   &  5.9  \\ \hline
		      LSP        &   FaNCL    & ${\bf 0.793} \pm 0.002$ &  13  & 1.9  & ${\bf 0.650} \pm 0.001$ &  9   &  6.1  \\ \hline
		      TNN        &   FaNCL    &    $0.798 \pm 0.001$    &  17  & 3.3  &    $0.655 \pm 0.002$    &  8   &  6.2  \\ \hline
	\end{tabular}
	\vspace{-5px}
	\label{tab:yahoo}
\end{table*}

%%%%%%%%%%%%%%%%%%%%%%%%%%%%%%%%%%%%%%%%%%%%%%%%

\subsection{Robust Principal Component Analysis}
\label{sec:expt2}

%%%%%%%%%%%%%%%%%%%%%%%%%%%%%%%%%%%%%%%%%%%%%%%%

\subsubsection{Synthetic Data}
\label{sec:synrpca}

\begin{table*}[ht]
	\centering
	\caption{RPCA performance of the various methods on synthetic data. The standard deviations of
		NMSE are all smaller than $0.0002$ and so not reported.
		CPU time is in seconds.}
	\vspace{-10px}
	\begin{tabular}{cc|ccc|ccc|ccc|ccc}
		\hline
		                &       &  \multicolumn{3}{c|}{$m=500$}   &  \multicolumn{3}{c|}{$m=1000$}  &  \multicolumn{3}{c|}{$m=1500$}  &  \multicolumn{3}{c}{$m=2000$}   \\
		                &       &     NMSE     & rank &   time    &     NMSE     & rank &   time    &     NMSE     & rank &   time    &     NMSE     & rank &   time    \\ \hline
		 nuclear norm   &  APG  &    $0.46$    &  5   &    1.5    &    $0.30$    &  10  &    9.7    &    $0.25$    &  15  &   33.9    &    $0.18$    &  20  &   94.7    \\ \hline\hline
		capped-$\ell_1$ &  GPG  & ${\bf 0.36}$ &  5   &    0.9    & ${\bf 0.25}$ &  10  &    6.7    & ${\bf 0.21}$ &  15  &   18.7    & ${\bf 0.15}$ &  20  &   60.4    \\ \cline{2-14}
		                & FaNCL & ${\bf 0.36}$ &  5   & {\bf 0.2} & ${\bf 0.25}$ &  10  & {\bf 1.4} & ${\bf 0.21}$ &  15  & {\bf 2.7} & ${\bf 0.15}$ &  20  & {\bf 6.5} \\ \hline
		      LSP       &  GPG  & ${\bf 0.36}$ &  5   &    2.7    & ${\bf 0.25}$ &  10  &   18.5    & ${\bf 0.21}$ &  15  &   111.2   & ${\bf 0.15}$ &  20  &   250.2   \\ \cline{2-14}
		                & FaNCL & ${\bf 0.36}$ &  5   & {\bf 0.4} & ${\bf 0.25}$ &  10  & {\bf 1.8} & ${\bf 0.21}$ &  15  & {\bf 3.9} & ${\bf 0.15}$ &  20  & {\bf 7.1} \\ \hline
		      TNN       &  GPG  & ${\bf 0.36}$ &  5   &    0.8    & ${\bf 0.25}$ &  10  &    6.0    & ${\bf 0.21}$ &  15  &   23.1    & ${\bf 0.15}$ &  20  &   51.4    \\ \cline{2-14}
		                & FaNCL & ${\bf 0.36}$ &  5   & {\bf 0.2} & ${\bf 0.25}$ &  10  & {\bf 1.2} & ${\bf 0.21}$ &  15  & {\bf 2.9} & ${\bf 0.15}$ &  20  & {\bf 5.8} \\ \hline
	\end{tabular}
	\label{tab:synperformance}
\end{table*}

In this section, 
we first perform experiments on a synthetic data set.  The observed $m \times m$ matrix is
generated as $O = U V + \tilde{Y} + G$, where
elements  of
$U \in \R^{m \times k}, V \in \R^{k \times m}$ (with $k = 0.01m$) are 
sampled i.i.d. from $\mathcal{N}(0, 1)$,
and elements of $G$ are sampled from $\mathcal{N}(0, 0.1)$.  Matrix $\tilde{Y}$ is sparse, with
$1\%$ of its elements randomly set to $5\| U V \|_{\infty}$  or $-5\| U V \|_{\infty}$ with equal probabilities. 
The sparsity regularizer
is 
the standard $\ell_1$, while different convex/nonconvex
low-rank regularizers are used.

For performance evaluation, 
we use
(i) NMSE 
$=\|(X+Y)-(UV+\tilde{Y})\|_F/\| UV+\tilde{Y}\|_F$, where
$X$ and $Y$
are the recovered low-rank and sparse components, respectively in (\ref{eq:rpca});
(ii) accuracy on locating the sparse support of $\tilde{Y}$ (i.e.,
percentage of entries that both $\tilde{Y}_{ij}$ and $Y_{ij}$ are nonzero or zero together);
and (iii) 
the recovered rank.
We vary $m$ in $\{500, 1000, 1500, 2000\}$.
Each experiment is repeated five times.

Note that IRNN  and the active subspace selection method cannot be used here.
Their objectives are of the form ``smooth function plus low-rank
regularizer", while RPCA  has a nonsmooth $\ell_1$ regularizer besides its 
low-rank
regularizer.
Similarly, Soft-Impute is
for matrix completion only.

Results are shown in Table~\ref{tab:synperformance}.
The accuracy on locating the sparse support
are always 100\% 
for all methods, and thus  are
not shown.
As can be seen, while both convex and nonconvex regularizers can perfectly recover the matrix rank and
sparse locations,
the nonconvex regularizers have lower NMSE's.
Moreover, as in matrix completion, FaNCL is again much faster. The larger
the matrix, the higher is 
the speedup.

%\begin{table}[ht]
%	\centering
%	\renewcommand{\arraystretch}{1.2}
%	\caption{Time comparison for various RPCA solvers on the synthetic data, with different matrix sizes $m$.}
%	\begin{tabular}{cc|c|c|c|c}
%		\hline
%		                   &       &   $m=500$    &     1000     &     1500     &     2000     \\ \hline
%		   nuclear norm    &  APG  &     1.5      &     9.7      &     33.9     &     94.7     \\ \hline
%		capped    $\ell_1$ &  GPG  &     0.9      &     6.7      &     18.7     &     60.4     \\ \cline{2-6}
%		                   & FaNCL & \textbf{0.2} & \textbf{1.4} & \textbf{2.7} & \textbf{6.5} \\ \hline
%		       LSP         &  GPG  &     2.7      &     18.5     &    111.2     &    250.2     \\ \cline{2-6}
%		                   & FaNCL & \textbf{0.4} & \textbf{1.8} & \textbf{3.9} & \textbf{7.1} \\ \hline
%		       TNN         &  GPG  &     0.8      &     6.0      &     23.1     &     51.4     \\ \cline{2-6}
%		                   & FaNCL & \textbf{0.2} & \textbf{1.2} & \textbf{2.9} & \textbf{5.8} \\ \hline
%	\end{tabular}
%	\label{tab:syntime}
%\end{table}

%%%%%%%%%%%%%%%%%%%%%%%%%%%%%%%%%%%%%%%%%%%%%%%%

\subsubsection{Background Removal on Videos}

\begin{figure}[H]
	\centering
	\subfigure[{\em bootstrap}.]
	{\includegraphics[width = 0.23\columnwidth]{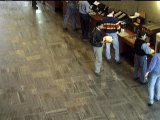}}
	\subfigure[{\em campus}.]
	{\includegraphics[width = 0.23\columnwidth]{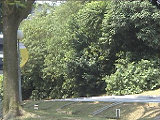}} 
	\subfigure[{\em escalator}.]
	{\includegraphics[width = 0.23\columnwidth]{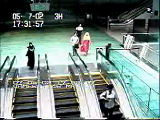}} 
	\subfigure[{\em hall}.]
	{\includegraphics[width = 0.23\columnwidth]{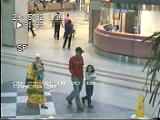}}
	\vspace{-.1in}
	\caption{Example image frames in the videos.}
	\label{fig:exampleImages}
	\vspace{-5px}
\end{figure}

In this section, we 
use RPCA
to perform 
video denoising on
background removal of 
corrupted videos.
Four benchmark 
videos\footnote{\url{http://perception.i2r.a-star.edu.sg/bk_model/bk_index.html}}  
in \cite{candes2011robust,sun2013robust}
are used
(Table~\ref{tab:sumVideo}), and example image frames are shown in
Figure~\ref{fig:exampleImages}.
As discussed in \cite{candes2011robust},
the stable image background can be treated as low-rank, while the 
foreground moving objects contribute to the 
sparse component.

\begin{table}[ht]
	\centering
	\vspace{-5px}
	\caption{Videos used in the experiment.}
	\vspace{-5px}
	\begin{tabular}{c | c | c | c | c}
		\hline
		& {\em bootstrap} & {\em campus} & {\em escalator} & {\em hall} \\ \hline
		\#pixels / frame & 19,200          & 20,480       & 20,800          & 25,344     \\ \hline
		total \#frames  & 9,165           & 4,317        & 10,251          & 10,752     \\ \hline
	\end{tabular}
	\vspace{-10px}
	\label{tab:sumVideo}
\end{table}

\begin{figure*}[t]
	\centering
	\subfigure[original.]
	{\includegraphics[width=0.38\columnwidth]{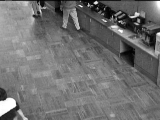}}
	\subfigure[nuclear norm.]
	{\includegraphics[width=0.38\columnwidth]{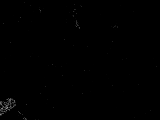}}
	\subfigure[capped-$\ell_1$.]
	{\includegraphics[width=0.38\columnwidth]{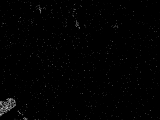}}
	\subfigure[LSP.]
	{\includegraphics[width=0.38\columnwidth]{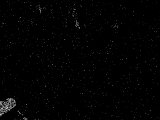}}
	\subfigure[TNN.]
	{\includegraphics[width=0.38\columnwidth]{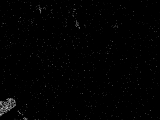}}
	\vspace{-.1in}
	\caption{Example foreground images in {\em bootstrap}, as recovered by using various low-rank
		regularizers.}
	\label{fig:foreground}
	\vspace{-10px}
\end{figure*}

Each image frame 
is reshaped as a column vector, and all frames are then stacked together to form a matrix.  
The pixel values are normalized to $[0, 1]$, and
Gaussian noise from $\mathcal{N}(0, 0.15)$ is added.
The experiment is repeated five times.

For performance evaluation, we use the commonly used peak signal-to-noise ratio 
\cite{dabov2007image}: PSNR $= - 10 \log_{10}(\text{MSE})$,
where $\text{MSE} = \frac{1}{m n} \sum_{i=1}^m \sum_{j=1}^n \left( X_{ij} - O_{ij} \right)^2$,
$X \in \R^{m \times n}$ is the recovered video,
and 
$O \in \R^{m \times n}$ 
is the ground-truth.

\begin{figure}[ht]
	\centering
	\vspace{-5px}
	\includegraphics[width = 0.3\textwidth]{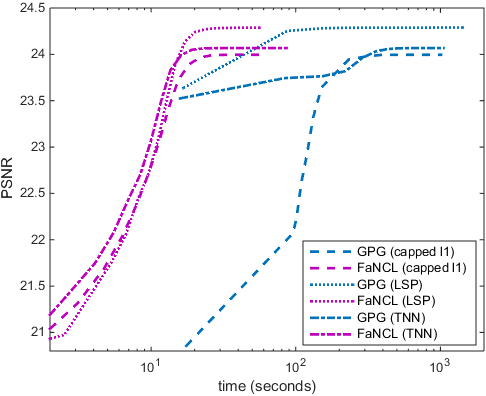}
	\vspace{-.1in}
	\caption{PSNR vs CPU time on the {\em bootstrap} data set.}
	\label{fig:bootstrap}
	\vspace{-15px}
\end{figure}

Results are shown in Table~\ref{tab:videoPSNR}.
As can be seen, the nonconvex regularizers lead to 
better  PSNR's than the convex nuclear norm. 
Moreover, FaNCL is more than $10$ times faster than GPG.
Figure~\ref{fig:foreground}
shows an example of the recovered foreground
in the {\em bootstrap} video.
As can been seen, the nonconvex regularizers can better separate foreground from background.
Figure~\ref{fig:bootstrap}
shows
the PSNR vs time on {\em bootstrap}.
Again, FaNCL converges much faster than others. 

\begin{table*}[ht]
	\centering
	\vspace{-10px}
	\caption{PSNR (in dB) and CPU time (in seconds) on the video background removal experiment. 
		For comparison, the PSNRs for all the input videos are 16.47dB.}
	\vspace{-10px}
	\begin{tabular}{cc|cc|cc|cc|cc}
		\hline
		\multicolumn{2}{c|}{}  &  \multicolumn{2}{c|}{{\em bootstrap}}   &    \multicolumn{2}{c|}{{\em campus}}    &  \multicolumn{2}{c|}{{\em escalator}}   &     \multicolumn{2}{c}{{\em hall}}      \\ \cline{3-10}
		\multicolumn{2}{c|}{}  &            PSNR            &    time    &            PSNR            &    time    &            PSNR            &    time    &            PSNR            &    time    \\ \hline
		nuclear norm   &  APG  &    $23.01 \!\pm\! 0.03$    &   688.4    &    $22.90 \!\pm\! 0.02$    &   102.6    &    $23.56 \!\pm\! 0.01$    &   942.5    &    $23.62 \!\pm\! 0.01$    &   437.7    \\ \hline\hline
		capped-$\ell_1$ &  GPG  &    $24.00 \!\pm\! 0.03$    &   1009.3   &    $23.14 \!\pm\! 0.02$    &    90.6    & ${\bf 24.33} \!\pm\! 0.02$ &   1571.2   &    $24.95 \!\pm\! 0.02$    &   620.0    \\ \cline{2-10}
		& FaNCL &    $24.00 \!\pm\! 0.03$    & {\bf 60.4} &    $23.14 \!\pm\! 0.02$    & {\bf 12.4} & ${\bf 24.33} \!\pm\! 0.02$ & {\bf 68.3} &    $24.95 \!\pm\! 0.02$    & {\bf 34.7} \\ \hline
		LSP       &  GPG  & ${\bf 24.29} \!\pm\! 0.03$ &   1420.2   & ${\bf 23.96} \!\pm\! 0.02$ &    88.9    &    $24.13 \!\pm\! 0.01$    &   1523.1   & ${\bf 25.08} \!\pm\! 0.01$ &   683.9    \\ \cline{2-10}
		& FaNCL & ${\bf 24.29} \!\pm\! 0.03$ & {\bf 56.0} & ${\bf 23.96} \!\pm\! 0.02$ & {\bf 17.4} &    $24.13 \!\pm\! 0.01$    & {\bf 54.5} & ${\bf 25.08} \!\pm\! 0.01$ & {\bf 35.8} \\ \hline
		TNN       &  GPG  &    $24.06 \!\pm\! 0.03$    &   1047.5   &    $23.11 \!\pm\! 0.02$    &   130.3    &    $24.29 \!\pm\! 0.01$    &   1857.7   &    $24.98 \!\pm\! 0.02$    &   626.2    \\ \cline{2-10}
		& FaNCL &    $24.06 \!\pm\! 0.03$    & {\bf 86.3} &    $23.11 \!\pm\! 0.02$    & {\bf 12.6} &    $24.29 \!\pm\! 0.01$    & {\bf 69.6} &    $24.98 \!\pm\! 0.02$    & {\bf 37.4} \\ \hline
	\end{tabular}
	\label{tab:videoPSNR}
	\vspace{-15px}
\end{table*}

%-----------------------------------------------------------------------------------------
\section{Conclusion}

In this paper, we considered the challenging problem of nonconvex low-rank matrix optimization.
The key observations are that for  the popular low-rank regularizers, 
the singular values 
obtained from the proximal operator
can be automatically thresholded, and also that the proximal operator can be
computed on a smaller matrix. For matrix completion, extra speedup can be achieved
by exploiting the ``sparse plus low-rank'' structure of the matrix estimate in each iteration.
The resultant algorithm is guaranteed to converge to a critical point of the nonconvex optimization problem.
Extensive experiments on matrix completion and RPCA 
show that the proposed algorithm is much faster than the state-of-art convex and nonconvex
low-rank solvers.
It also demonstrates that nonconvex low-rank regularizers outperform the 
convex nuclear norm regularizer in terms of recovery accuracy and the rank obtained.

\section*{Acknowledgment}

This research was supported in part by
the Research Grants Council of the Hong Kong Special Administrative Region
(Grant 614513).

%-----------------------------------------------------------------------------------------
% bibliography
%-----------------------------------------------------------------------------------------
{
\footnotesize
\bibliographystyle{IEEEtran}
\bibliography{bib}
}

\cleardoublepage
\appendix

%%%%%%%%%%%%%%%%%%%%%%%%%%%%%%%%%%%%%%%%%%%%%%%%

\subsection{Proof of Proposition~\ref{pr:proxReduce}}
\label{app:proRed}

This Proposition appears in \cite{lu2015generalized,gu2014weighted},
for completeness, we also present a proof here.
First, note that for a matrix $X \in \R^{m \times n}$ and any orthogonal projection matrices 
$P \in \R^{m \times m}$ and 
$Q \in \R^{n \times n}$ (where
$P^{\top}P = I, 
Q^{\top}Q = I$),
$X$ has the same singular values with $P^{\top} X Q$.
Therefore, Assumption A3 implies $r$ is invariant to orthogonal projection, i.e.,
\begin{align*}
r(X) = r(P^{\top} X Q)
\end{align*}
Then, we introduce the following Proposition in \cite{boyd2004convex}.
\begin{prop} \label{pr:maxSv}
	Let $(u^*,v^*) = \arg \max_{u, v} (u^{\top} X v)\;:\; \SN{u} = \SN{v} = 1$.  
	Then, $u^*$ (resp.
	$v^*$)  is the left (resp. right) singular vector of $X$, and
	the optimal objective value is $\sigma_1$, the largest singular value of $X$.
\end{prop}

Let the SVD of $X$ be $P \Sigma_X Q^{\top}$.
Since, $\|\cdot\|_F$
and $r(\cdot)$ are invariant
to orthogonal projections,
\begin{eqnarray}
\mathcal{P}^{\mu}_{r(\cdot)}(Z) & =  & \min_X \frac{1}{2} \FN{X - Z}^2 + \mu r(X) 
\label{eq:app2} \\
& = &  \min_{P, \Sigma_X, Q} \frac{1}{2}\FN{P\Sigma_X Q^{\top} - Z}^2 
+ \mu r(P\Sigma_X Q^{\top})
\notag \\
& = &  \min_{P, \Sigma_X, Q} \frac{1}{2}\FN{\Sigma_X - P^{\top}ZQ}^2 
+ \mu r(\Sigma_X)
\notag \\
&= & \min_{\Sigma_X}\frac{1}{2}\Tr{\Sigma_X^2 + Z^{\top}Z} + \mu r(\Sigma_X) 
\notag \\
&& - \max_{P, Q} \Tr{\Sigma_X P^{\top}ZQ}.
\notag
\end{eqnarray}
Let $p_i$ (resp. $q_i$) be the $i$th column of $P$ (resp. $Q$). 
We have
\[ \max_{P, Q} \Tr{\Sigma_X P^{\top}ZQ}
= \sum_{i=1}^n \left[\sigma_X \right]_i 
\max_{p_i, q_i}p_i^{\top} Z q_i. \]
Recall that the SVD of $Z$ is $U \Sigma V^{\top}$.
Using Proposition~\ref{pr:maxSv}, 
$\sigma_1 = \max_{p_1, q_1}p_1^{\top} Z q_1$, and
$p_1 = u_1, q_1 = v_1$ where $u_i$ (resp. $v_i$) is the $i$th column of $U$ (resp. $V$). 
Since $p_i \neq p_j, q_i \neq q_j$ if $i \neq j$, again by Proposition~\ref{pr:maxSv},
we have $p_2 = u_2, q_2 = v_2$, and so on. Hence,
$P = U, Q = V$ and $\Sigma = P^{\top}ZQ$, then 
(\ref{eq:app2}) can then be rewritten as:
\begin{eqnarray}
\lefteqn{\min_{\Sigma_X} \frac{1}{2} \FN{\Sigma_X - \Sigma}^2 + \mu r(\Sigma_X)}
\label{eq:tmp}\\
& = & \min_{y \ge 0} \frac{1}{2} \sum_{i = 1}^n (y_i - \sigma_i)^2 + \mu r(\Diag(y)),
\nonumber
\end{eqnarray}
which leads to (\ref{eq:proRed}).

%%%%%%%%%%%%%%%%%%%%%%%%%%%%%%%%%%%%%%%%%%%%%%%%
\subsection{Proof of Proposition~\ref{pr:proxSolution}}
\label{app:proxSolution}

\subsubsection{}
First, we show $y^*_i \le \sigma_i$.
By assumption A3, (\ref{eq:tmp}) 
(or, equivalently, (\ref{eq:proRed}) with Proposition~\ref{pr:proxReduce})
can be rewritten as
\begin{eqnarray*}
	\lefteqn{\min_{\Sigma_X} \frac{1}{2}\FN{\Sigma_X - \Sigma}^2 + \mu \sum_{i = 1}^n
		\hat{r}\left(\left[\sigma_X\right]_i\right)} \\
	& = & \sum_{i = 1}^n \min_{y_i \ge 0} \frac{1}{2} (y_i - \sigma_i)^2 + \mu \hat{r}(y_i) , 
\end{eqnarray*}
If the optimal $y_i^*$ is achieved at the boundary (i.e., $y_i^* = 0$), 
obviously $y_i^* \le \sigma_i$.
Otherwise, 
\begin{equation} \label{eq:app3}
0 \in y_i^* - \sigma_i + \mu\partial \hat{r}(y_i^*).
\end{equation} 
Since $\hat{r}(x)$ is non-decreasing on $x \ge 0$,
its super-gradient $\partial \hat{r}(\cdot)$ is non-negative on $\R^+$, and so $y_i^* \le \sigma_i$.

\subsubsection{}
Now, consider an $(i,j)$ pair such that 
$\sigma_j \geq \sigma_i$.
Assume that $y^*_i\geq 0$ and $y^*_j\geq 0$. 
From (\ref{eq:app3}), we have
\begin{equation*}
0 \in y_i^* - \sigma_i + \mu\partial \hat{r}(y_i^*) \text{ and }
0 \in y_j^* - \sigma_j + \mu\partial \hat{r}(y_j^*).
\end{equation*} 

Again from assumption A3, since $\hat{r}(x)$ is concave and non-decreasing on $x \ge 0$,
its super-gradient $\partial \hat{r}(\cdot)$ is thus also non-increasing.
If $y_i^*$ is not achieved at the boundary (i.e., is locally optimal), consider $\sigma_j >
\sigma_i > 0$. To ensure that (\ref{eq:app3}) holds, we can ether (i) $y_j^* > y_i^*$, 
and thus $\partial \hat{r}(y_j^*) \le \partial \hat{r}(y_i^*)$; or (ii) $y_j^* < y_i^*$, and thus $\partial
\hat{r}(y_j^*) \ge \partial \hat{r}(y_i^*)$.  
However, $\partial \hat{r}(\cdot)$ is non-negative, 
and thus lower-bounded.
Hence, there always exists $y_j^* > y_i^*$ to ensure (\ref{eq:app3}).
If multiple solutions exist, we
take the largest one. So, we must have $y_j^* > y_i^*$.  

\subsubsection{}
Thus, the smaller the $\sigma_i$, the smaller is $y_i^*$ and $y_i^* \le \sigma_i$. Since
$\hat{r}(y_i^*)$ is non-increasing on $\R^+$, $\mu\hat{r}(y_i^*)$ will not become
smaller. Thus, there must exists $\gamma$ such that once $\sigma_i \le \gamma$,
(\ref{eq:app3}) no longer holds, and $y_i^*$ is not locally optimal and lies on
the boundary (i.e., $y_i^* = 0$).

\subsection{Proof of Proposition~\ref{pr:approGSVT}}
\label{app:approGSVT}

Since $\Span{U_{\hat{k}}} \subseteq \Span{Q}$ and $Q$ is orthogonal, it can be written as
$Q = [U_=; U_{\bot}] R$
where $\Span{U_=} = \Span{U_{\hat{k}}}$, $U_{\bot}^{\top} U_{\hat{k}} = 0$
and $R$ is a rotation matrix ($R R^{\top} = R^{\top} R = I$).
Thus, $Q^{\top} Z^t = R^{\top} [U_{=}; U_{\bot}]^{\top} Z^t$
and its rank-$\hat{k}$ SVD is 
$R^{\top} [U_{=}; 0]^{\top} U_{\hat{k}} \Sigma V_{\hat{k}}^{\top}$.
Using Proposition~\ref{pr:proxReduce}, 
\begin{align*}
\Prox{\mu}{r(\cdot)}{Q^{\top} Z^t} = 
R^{\top} 
\begin{bmatrix}
U_{=}^{\top} \\
0
\end{bmatrix} U_{\hat{k}} \hat{\Sigma} V^{\top}_{\hat{k}},
\end{align*}
where $\hat{\Sigma} = \Diag(y_1^*, \dots, y_{\hat{k}}^*)$ is the optimal solution in
(\ref{eq:proRed}). 
Then note that,
\begin{align*}
Q \Prox{\mu}{r(\cdot)}{Q^{\top} Z^t} 
& = [U_{=}; U_{\bot}] R R^{\top}
\begin{bmatrix}
U_{=}^{\top} \\
0
\end{bmatrix}   U_{\hat{k}} \hat{\Sigma} V_{\hat{k}}^{\top} \\
& = [U_{=}; U_{\bot}]
\begin{bmatrix}
U_{=}^{\top} \\
0
\end{bmatrix}   U_{\hat{k}} \hat{\Sigma} V_{\hat{k}}^{\top} \\
& = U_{=} U_{=}^{\top} U_{\hat{k}} \hat{\Sigma} V^{\top}_{\hat{k}}.
\end{align*}

Since $\Span{U_{\hat{k}}} = \Span{U_=}$, so $U_= U_=^{\top}= U_{\hat{k}}
U_{\hat{k}}^{\top}$, and 
%note that $(U_{\hat{k}}^{\top} U_{\hat{k}}) = I$, so
\begin{align*}
U_{=} U_{=}^{\top} U_{\hat{k}} \hat{\Sigma} V_{\hat{k}}^{\top}
= U_{\hat{k}} (U_{\hat{k}}^{\top} U_{\hat{k}}) \hat{\Sigma} V_{\hat{k}}^{\top}
= U_{\hat{k}} \hat{\Sigma} V_{\hat{k}}^{\top},
\end{align*}
which is $\Prox{\mu}{r(\cdot)}{Z^t}$.

%%%%%%%%%%%%%%%%%%%%%%%%%%%%%%%%%%%%%%%%%%%%%%%%

\subsection{Proof of Proposition~\ref{pr:bound}}
\label{app:convergence}

Since (\ref{eq:decrease}) holds, we have
\begin{align*}
F(X^{t+1}) \le F(X^{t}) - c_1 \FN{X^{t+1}-X^{t}}.
\end{align*}

Sum it from $0$ to $T$,  we have
\begin{align}
F(X^0) - F(X^{T + 1}) \ge c_1 \sum_{t = 1}^{T}\FN{X^{t+1}-X^t}^2.
\label{eq:temp1}
\end{align}

By assumption A2, $F(X)$ is bounded below.
Thus, as $T \rightarrow + \infty$, there exists a finite constant $\alpha$ such that
\begin{align}
\alpha = \sum_{t = 1}^{+ \infty} \FN{X^{t+1}-X^t}^2.
\label{eq:app8}
\end{align}

Hence, we must have
$\lim\limits_{t \rightarrow \infty}\FN{X^{t+1} - X^t}^2 = 0$, and
thus $\{X^t\}$ converges to a limit point $X^*$. 

\subsection{Proof of Theorem~\ref{the:convergence}}

Next, we show that $X^*$ is a critical point of (\ref{eq:problem}).
First, as in \cite{gong2013general}, it is easy to see that $r(\cdot)$ here can also be
decomposed as the difference 
of two convex functions $r_1(X)$ and $r_2(X)$ i.e.,
$r(X) = r_1(X) - r_2(X)$.
Consider the optimal conditions in proximal step, we have
\begin{align}
0 \; \in \; & \nabla \f{X^t} + X^{t + 1} - X^t
\label{eq:app9} \\
& + \lambda \partial r_1(X^{t+ 1}) - \lambda \partial r_2(X^{t+ 1}).
\notag
\end{align}

For limit point $X^{t + 1} = X^t = X^*$, so $X^{t + 1} - X^t = 0$ and vanish. Thus,
\begin{align*}
0 \in \nabla \f{X^*} + \lambda \partial r_1(X^*) - \lambda \partial r_2(X^*),
\end{align*}
and
$X^*$ is a critical point of (\ref{eq:problem}).

\subsection{Proof of Corollary~\ref{the:conv:rate}}

In Theorem~\ref{the:convergence},
we have shown Algorithm~\ref{alg:FaNCL} can converge to a critical point of \eqref{eq:problem}.
Then, from \eqref{eq:temp1}, 
rearrange items we will have
\begin{align*}
\min_{t = 1, \cdots, T} \FN{X^{t+1}-X^t}^2
& \le
\frac{1}{T}\sum_{t = 1}^{T}\FN{X^{t+1}-X^t}^2 \\
& \le
\frac{1}{c_1 T}\left[ F(X^1) - F(X^{T + 1})  \right],
\end{align*}
which proves the Corollary.

%%%%%%%%%%%%%%%%%%%%%%%%%%%%%%%%%%%%%%%%%%%%%%%%

\subsection{Proof of Proposition~\ref{pr:checkobject}}
\label{app:checkobject}

By definition of the Frobenoius norm, 
we only need to show that $P^{\top} X Q$ has the same singular values as $X$. 
Since $U \subseteq P$, we partition $P$ as $P = [U_{=}, U_{\bot}] R$,
where $\Span{U_{=}} = \Span{U}$, $U_{\bot}$ is orthogonal to $U$ (i.e., $U_=^{\top} U_{\bot}
= 0$),
and $R$ is a rotation matrix.
Then,
\begin{align*}
P^{\top} X 
= R^{\top} \left[ U_{=}, U_{\bot} \right]^{\top} U \Sigma V^{\top}
= R^{\top} [U_{=}, 0]^{\top} U \Sigma V^{\top}.
\end{align*}
Since $\Span{U_{=}} = \Span{U}$, we have $U_{=} U_=^{\top} = U U^{\top}$. 
Let $\hat{U} = R^{\top} [U_{=}, 0]^{\top} U$. Then,
\begin{align*}
\hat{U}^{\top} \hat{U} 
& = U^{\top} [U_{=}, 0] R R^{\top} [U_{=}, 0]^{\top} U \\
& = U^{\top} U_{=} U_{=}^{\top} U 
  = (U^{\top} U) (U^{\top} U) = I.
\end{align*}
Hence, $\hat{U} \Sigma V^{\top}$ is the reduced SVD of $P^{\top}X$. 
Similarly, for $Q$, we obtain that $\Sigma$ is also the singular values of $P^{\top} X Q$.

%%%%%%%%%%%%%%%%%%%%%%%%%%%%%%%%%%%%%%%%%%%%%%%%

\subsection{Proof of Theorem~\ref{the:conv:altprox}}
\label{app:convergence:rpca}

Here, we prove the case for two blocks of parameters \eqref{eq:rpca} as an example.
Extension to multiple block is easily obtained. 

\subsubsection{}
Let
$\Delta_{X^t}^2 = \FN{X^{t+1} - X^{t}}$, and
$\Delta_{Y^t}^2 = \FN{Y^{t+1} - Y^{t}}$.
When sufficient decrease holds for both $X$ and $Y$, we have
\begin{eqnarray*}
	F(X^{t+1},Y^{t+1}) 
	\!\!\!\! & \le & \!\!\!\! F(X^{t+1},Y^t) - c_1 \Delta_{Y^t}^2 \notag\\
	\!\!\!\! & \le & \!\!\!\! F(X^t,Y^t) - c_1 \Delta_{X^t}^2 - c_1 \Delta_{Y^t}^2
\end{eqnarray*}

Summarize above from $t = 0$ to $T$, we get
\begin{eqnarray}
F(X^1, Y^1) - F(X^{T + 1}, Y^{T + 1}) 
\ge c_1 \sum_{t = 0}^T \left( \Delta_{X^t}^2 + \Delta_{Y^t}^2 \right).
\label{eq:temp2}
\end{eqnarray}

Since $F(X, Y)$ is bounded below, L.H.S above is a finite positive constant.
Same as (\ref{eq:app8}):
\begin{align*}
\lim\limits_{t \rightarrow \infty}\FN{X^{t+1} - X^t} = 0, \quad
\lim\limits_{t \rightarrow \infty}\FN{Y^{t+1} - Y^t} = 0.
\end{align*}

Thus, $\sum_{t = 0}^T \left( \Delta_{X^t}^2 + \Delta_{Y^t}^2 \right) \le + \infty$.

\subsubsection{}
From the optimal conditions of the proximal step, similar to (\ref{eq:app9}), we have
\begin{align*}
0 & \in \nabla_Y f( X^*, Y^*) + \beta \partial \|Y^*\|_1, \\
0 & \in \nabla_X f( X^*, Y^*) + \lambda \partial r_1 (X^*) - \lambda \partial r_2 (X^*).
\end{align*}
Thus, $(X^*, Y^*)$ is a critical point of (\ref{eq:rpca}).

\subsubsection{}
Finally, using same technique at proof of Corollary~\ref{the:conv:rate} and \eqref{eq:temp2},
it is easy to obtain
\begin{align*}
\min_{t = 1,\cdots,T} \left( \Delta_{X^t}^2 + \Delta_{Y^t}^2 \right)
\le \frac{1}{c_1 T}\left[ F(X^1, Y^1) - F(X^{T + 1}, Y^{T + 1})  \right].
\end{align*}

%%%%%%%%%%%%%%%%%%%%%%%%%%%%%%%%%%%%%%%%%%%%%%%%%%%%%%%%%%%%%%%%%%

\subsection{Solution of GSVT and details of Corollary~\ref{cor:proxSolution}}
\label{app:solGSVT}

To simplify notations, 
%on one dimension case, 
we write $y_i$ as $y$, and $\sigma_i$ as $\sigma$.
Our focus here is $\gamma$ and is derived based on GIST \cite{gong2013general}.
For LSP, MCP and SCAD, 
the relationship between different stationary points is ignored in \cite{gong2013general}, 
thus their solutions are not necessarily the local optimal.

%Closed-form solution for (\ref{eq:app12}) has also been provided by \cite{gong2013general}.
%However, they ignore the relationship among $h(0)$, $h(\hat{y}_1)$ and $h(\hat{y}_2)$. 
%Hence, their "solution" is only a stationary point, but not necessarily the optimal solution.

%%%%%%%%%%%%%%%%%%%%%%%%%%%%%%%%%%%%%%%%%%%%%%%%

\subsubsection{$\ell_1$-regularizer}

The closed-form solution is at Lemma~\ref{lem:svt}, 
as it can be seen that $\gamma = \mu$ for the nuclear norm.
%\footnote{$\surd$ dont know what u're trying to say. if "Close-form solution is given at Lemma~\ref{lem:svt}", then it's done. what's the purpose of the other sentences???}

%%%%%%%%%%%%%%%%%%%%%%%%%%%%%%%%%%%%%%%%%%%%%%%%

\subsubsection{LSP}

For LSP, (\ref{eq:proRed}) becomes
\begin{align*}
\min_y h(y) \equiv \frac{1}{2}(y-\sigma)^2 + \mu \log\left(1 + \frac{y}{\theta}\right).
\end{align*}
If $\sigma = 0$, obviously $y^* = 0$. 
So we only need to consider $\sigma > 0$. Now,
\begin{align*}
\nabla h(y) = y - \sigma + \frac{\mu}{y + \theta}.
\end{align*}
Since $\theta, y > 0$, 
\begin{align}
(\theta + y)\nabla h(y) 
& = (y + \theta)(y - \sigma) + \mu. 
\notag \\
& = y^2 - (\sigma - \theta)y + \mu - \theta\sigma.
\label{eq:app10}
\end{align}
%and is a simple quadratic problem with discriminant

\textbf{Case 1: 
	$\Delta \equiv (\sigma + \theta)^2 - 4 \mu \le 0$}: 
Then $\nabla h(y) \ge 0$ on $\R^+$, and 
thus $h(y)$ is non-deceasing on $y \ge 0$. If
$0 \le \sigma \le \min\left(0, -\theta + 2 \sqrt{\mu}\right)$,
we have $\arg\min_y h(y) = 0$.

\begin{figure}[H]
	\centering
	\vspace{-5px}
	\includegraphics[width = 120px]{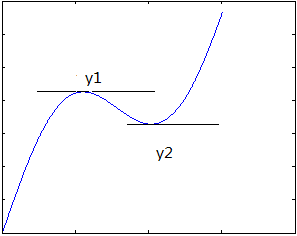}
	\vspace{-5px}
	\caption{Illustration for Case 2.}
	\label{fig:app1}
	\vspace{-5px}
\end{figure}

\textbf{Case 2: $\Delta > 0$.} 
The square roots of 
$y^2 - (\sigma - \theta)y + \mu - \theta\sigma=0$
in 
(\ref{eq:app10}) are
\begin{eqnarray*}
	\hat{y}_1 & = & \frac{1}{2}\left(\sigma-\theta-\sqrt{(\sigma+\theta)^2-4\mu}\right), \\
	\hat{y}_2 & = & \frac{1}{2}\left(\sigma-\theta+\sqrt{(\sigma+\theta)^2-4\mu}\right).
\end{eqnarray*}
%which are two stationary points. 
Since $h(y)$ has two stationary points, it is of the form in Figure~\ref{fig:app1},
and $y^*$ depends only on 
\begin{align*}
h(0)   = \frac{1}{2} \sigma^2, \quad
h(\hat{y}_2) = h\left( \frac{1}{2}(\sigma-\theta + \sqrt{\Delta}) \right).
\end{align*}
Thus, if 
$h(0) < h(\hat{y}_2)$, $y^* = 0$. When $h(0) = h(\hat{y}_2)$, we take the largest one as
$y^* = \max \left(0, \hat{y}_2\right)$ (and thus the solution may not be unique).
Finally, 
when $h(0) > h(\hat{y}_2)$, we have $y^* = \hat{y}_2$.

However, obtaining $\gamma$ 
by directly comparing $h(0)$ and $h(\hat{y}_2)$ is complicated
and has no simple closed-form solution. 
Here, we take a simpler approach. Once $\hat{y}_2 \le 0$, we have $y^* = 0$. I.e., if
$\sigma \le \min \left( \theta, \frac{\mu}{\theta} \right)$,
we have $\hat{y}_2 \le 0$ and $y^* = 0$.

Finally, on combining both cases, the threshold for LSP can be obtained as
\begin{eqnarray*}
	\gamma & = & \max 
	\left\lbrace
	\min\left(0, -\theta + 2 \sqrt{\mu}\right),  \min\left(\frac{\mu}{\theta},\theta\right)
	\right\rbrace \\ 
	& = & \min\left(\frac{\mu}{\theta},\theta\right). 
\end{eqnarray*}
Using Proposition~\ref{pr:proxReduce}, the optimal solution is shown in Lemma~\ref{lem:proLog}.

\begin{lemma} \label{lem:proLog}
	When $r(\cdot)$ is the LSP,
	the optimal solution of 
	the corresponding proximal operator 
	is $\Prox{\mu}{r(\cdot)}{Z}= U\Diag(y_1^*,\dots,y_n^*)V^{\top}$
	where
	\begin{align*}
	y_i^* = 
	\begin{cases}
	0 & \sigma_i \le \min(0, -\theta + 2 \sqrt{\mu}), \\
	0 & \sigma_i > \min(0, -\theta + 2 \sqrt{\mu})   \text{\;and\;} h(0) < h(\hat{y}_2), \\
	\hat{y}_2 & \sigma_i > \min(0, -\theta + 2 \sqrt{\mu}) \text{\;and\;} h(0) \ge h(\hat{y}_2).
	\end{cases}
	\end{align*}
	depends only on $\sigma_i$,
	and $\hat{y}_2 = \frac{1}{2}(\sigma_i - \theta + \sqrt{(\sigma + \theta)^2 - 4 \mu})$.
\end{lemma}

\subsubsection{Capped $\ell_1$}

Problem (\ref{eq:proRed}) then becomes
\begin{align*}
h(y) \equiv \frac{1}{2}\left(y - \sigma\right)^2 + \mu \min\left(y, \theta\right).
\end{align*}
This can be written as 
\[ \arg\min h(y) =\left\{ \begin{array}{ll} 
\arg\min h_1(y) & y \le \theta \\
\arg\min h_2(y) & y >\theta 
\end{array} \right., \]
where
$h_1(y) = \frac{1}{2}(y-\sigma)^2 + \mu y$,     and 
$h_2(y) = \frac{1}{2}(y-\sigma)^2 + \mu \theta$.
The optimal cases for $h_1(y)$ are:
\begin{align*}
\begin{cases}
h_1(y^* = 0) = 0, 
& \sigma - \mu \le 0,\\
h_1(y^* = \sigma - \lambda) = -\frac{1}{2}\mu^2 + \mu \sigma, 
& 0 < \sigma - \mu < \theta,\\
h_1(y^* = \theta) = \frac{1}{2}\left(\theta-\sigma\right)^2 + \mu \sigma,   
& \theta \le \sigma - \mu, 
\end{cases}
\end{align*}
And those for $h_2(y)$ are:
\begin{align*}
\begin{cases}
h_2(y^* = \theta) = \frac{1}{2}\left(\theta-\sigma\right)^2 + \mu \theta,
& \sigma \le \theta, \\
h_2(y^* = \sigma) = \mu \theta,  
& \sigma > \theta.
\end{cases}
\end{align*}
Consider cases of $\theta \le \mu$ and $\theta > \mu$.
Taking the minimum over above functions, the optimal of $x^*$ is:
\begin{align*}
y^* =
\begin{cases}
(\sigma - \mu)_{+} & \sigma \le \theta + \frac{1}{2}\mu, \\
\theta            & \sigma >   \theta + \frac{1}{2}\mu.
\end{cases}
\end{align*}
Thus, for the capped-$\ell_1$ regularizer, $\gamma = \min ( \mu, \frac{1}{2}\lambda + \theta)$.
Combining with Proposition~\ref{pr:proxReduce}, we obtain the following:

\begin{lemma}
	\label{lem:proCAP}
	When $r(\cdot)$ is the	capped-$\ell_1$,
	the optimal solution of the corresponding proximal operator is
	$\Prox{\mu}{r(\cdot)}{Z} = U\Diag(y_1^*,\dots,y_n^*)V^{\top}$, where
	\begin{align*}
	y_i^* =
	\begin{cases}
	(\sigma_i - \mu)_{+} & \sigma_i \le \theta + \frac{1}{2}\mu \\
	\theta               & \sigma_i >   \theta + \frac{1}{2}\mu
	\end{cases},
	\end{align*}
	and depends on $\sigma_i$.
\end{lemma}

%%%%%%%%%%%%%%%%%%%%%%%%%%%%%%%%%%%%%%%%%%%%%%%%

\subsubsection{TNN}
For the TNN, it directly controls the number of singular values. However,
from Lemma~\ref{lem:proTNN}, it is to see that $\gamma = \min\left(\mu, \sigma_{\theta + 1}\right)$.
\begin{lemma}
	\cite{hu2013fast}
	\label{lem:proTNN}
	When $r(\cdot)$ is the TNN regularizer,
	the optimal solution of the proximal operator is
	\begin{align*}
	\Prox{\mu}{r(\cdot)}{Z}
	= U\left(\Sigma - \mu \tilde{I}_{\theta} \right)_{+}V^{\top},
	\end{align*}
	where $\tilde{I}_k$ is the square matrix with all zeros elements except 
	at positions $[\tilde{I}_{\theta}]_{ii} = 1$ for $i > \theta$.
\end{lemma}

%%%%%%%%%%%%%%%%%%%%%%%%%%%%%%%%%%%%%%%%%%%%%%%%

\subsubsection{MCP}

For MCP, again, $y^*$ for problem (\ref{eq:proRed}) becomes
\begin{align*}
y^* = 
\begin{cases}
\arg\min h_1(y) & 0 \le y \le  \theta \mu \\
\arg\min h_2(y) & y  >  \theta\mu
\end{cases}
\end{align*}
where $h_1(y)$ and $h_2(y)$ are defined as
\begin{align*}
h_1(y) & = \frac{1}{2}(1 - \frac{1}{\theta})y^2 - (\sigma - \mu) y + \frac{1}{2} \sigma^2  \\
h_2(y) & = \frac{1}{2}y^2 - \sigma y + \frac{1}{2}\sigma^2 + \frac{1}{2}\theta \mu^2   
\end{align*}

For $h_1(y)$, the optimal depends on $\theta$ as:
\begin{itemize}
	\item[(1).] If $\theta = 1$, then the optimal is:
	\begin{align*}
	y^* = 
	\begin{cases}
	0          & 0 \le \sigma \le \mu \\
	\mu & \sigma  >  \mu 
	\end{cases}
	\end{align*}
	
	\item[(2).] If $\theta > 1$, note that it is a quadratic function and the optimal depends on
	$\bar{y} = \frac{\theta (\sigma - \mu)}{\theta - 1}$. As a result:
	\begin{align*}
	y^* =
	\begin{cases}
	0 &  0 \le \sigma \le \mu  \\
	\bar{y} & \mu < \sigma < \theta\mu \\
	\theta \mu & \sigma \ge \theta\mu
	\end{cases}
	\end{align*}
	
	\item[(3).] If $0 < \theta < 1$, again it is a quadratic function, but the coefficient on
	the quadratic term is negative. Thus
	\begin{align*}
	y^* =
	\begin{cases}
	0 & 0 \le \sigma \le \frac{1}{2}\theta\mu + \frac{1}{2}\mu \\
	\theta\mu & \sigma >   \frac{1}{2}\theta\mu + \frac{1}{2}\mu
	\end{cases}
	\end{align*}
\end{itemize}
Then, for $h_2(y)$, it is simple:
\begin{align*}
y^* =
\begin{cases}
\theta\mu & 0 \le \sigma \le \theta\mu \\
\sigma & \sigma  >  \theta\mu 
\end{cases}
\end{align*}
Combine $h_1(y)$ and $h_2(y)$:
\begin{itemize}
	\item[(1).] If $\theta = 1$, then
	\begin{align}
	y^* =
	\begin{cases}
	0      & 0 \le \sigma \le \mu \\
	\sigma & \sigma > \mu
	\end{cases}
	\label{eq:app24}
	\end{align} 
	
	\item[(2).] If $\theta > 1$, then ($\bar{y} = \frac{\theta (\sigma - \mu)}{\theta - 1}$):
	\begin{eqnarray}
	y^* =
	\begin{cases}
	0 & 0 \le \sigma \le \mu \\
	\bar{y} & \mu <\sigma< \theta\mu \text{\;and\;} h_1(\bar{y}) \le h_2(\theta\mu) \\
	\theta\mu  & \mu <\sigma< \theta\mu \text{\;and\;} h_1(\bar{y}) > h_2(\theta\mu) \\
	\sigma & \sigma \ge \theta\mu 
	\end{cases}
	\label{eq:app25}
	\end{eqnarray}
	
	\item[(3).] If $0 < \theta < 1$, we need to compare $h_1(0)$ and $h_2(\sigma)$,
	then we have:
	\begin{align}
	y^* =
	\begin{cases}
	0      & 0 \le \sigma \le \sqrt{\theta}\mu \\
	\sigma & \sigma > \sqrt{\theta}\mu 
	\end{cases}
	\label{eq:app26}
	\end{align}
\end{itemize}
Thus, $\gamma$ for MCP is:
\begin{align*}
\gamma = 
\begin{cases}
\sqrt{\theta} \mu & 0 < \theta < 1 \\
\mu               & \theta \ge 1 
\end{cases}
\end{align*}
Using Proposition~\ref{pr:proxReduce}, the optimal solution is shown in Lemma~\ref{lem:proMCP}.
\begin{lemma} \label{lem:proMCP}
	When $r(\cdot)$ is the MCP, the 
	optimal solution of the corresponding proximal operator is:
	\begin{align*}
	\Prox{\mu}{r(\cdot)}{Z}
	= U\Diag(y_1^*,\dots,y_n^*)V^{\top},
	\end{align*}
	where $y_i^*$ depends on $\theta$ and $\sigma_i$, i.e. 
	if $\theta > 1$, then $y_i^*$ is given by (\ref{eq:app24});
	then if $\theta = 1$, then $y_i^*$ is given by (\ref{eq:app25});
	finally, if $0 < \theta < 1$, $y_i^*$ is given by (\ref{eq:app26}).
\end{lemma}

\subsubsection{SCAD}

Again, it can be written as $(\theta > 2)$
\begin{align*}
y^* = 
\begin{cases}
\arg\min h_1(y) & 0 \le y \le \mu \\
\arg\min h_2(y) & \mu < y \le \theta\mu \\
\arg\min h_3(y) & \theta\mu <  y  
\end{cases}
\end{align*}
where $h_1(y)$, $h_2(y)$ and $h_3(y)$ are defined as
\begin{align*}
h_1(y) & = \frac{1}{2}(y - \sigma)^2 + \mu y, \\
h_2(y) & = \frac{1}{2}(y - \sigma)^2 +  \frac{-y^2 + 2 \theta \mu y - \mu^2}{2(\theta - 1)}, \\
h_3(y) & = \frac{1}{2}(y - \sigma)^2 + \frac{(\theta + 1) \mu^2}{2}.
\end{align*}
Thus, we can get
\begin{itemize}
	\item[(1).] For $h_1(y)$, the optimal is
	\begin{align*}
	y^* =
	\begin{cases}
	0      & 0 \le \sigma \le \mu \\
	\sigma - \mu & \sigma >   \mu 
	\end{cases}
	\end{align*}
	
	\item[(2).] For $h_2(y)$, the optimal is
	\begin{align*}
	y^* =
	\begin{cases}
	2\mu & 0 \le \sigma \le 2\mu \\
	\frac{(\theta - 1)\sigma - \theta \mu}{\theta - 2} & 2\mu < \sigma < \theta \mu \\
	\theta \mu & \sigma \ge \theta \mu
	\end{cases}
	\end{align*}
	
	\item[(3).] For $h_3(y)$, the optimal is
	\begin{align*}
	y^* =
	\begin{cases}
	\theta \mu & 0 \le \sigma \le \theta\mu \\
	\sigma & \sigma >  \theta\mu 
	\end{cases}
	\end{align*}
\end{itemize}

To get $\gamma$, we need to compare $h_1(0)$, $h_2(\mu)$ and $h_3(\theta\mu)$, 
it is easy to verify $h_1(0)$ is smallest,
thus $\gamma = \mu$.
Finally, using Proposition~\ref{pr:proxReduce}, the optimal solution is shown in Lemma~\ref{lem:proSCAD}.
\begin{lemma} \label{lem:proSCAD}
	When $r(\cdot)$ is the SCAD, the 
	optimal solution of the corresponding proximal operator is
	$\Prox{\mu}{r(\cdot)}{Z}= U\Diag(y_1^*,\dots,y_n^*)V^{\top}$
	where
	\begin{align*}
	y_i^* = 
	\begin{cases}
	0 & 0 \le \sigma_i \le \mu \\
	\sigma_i - \mu & \mu < \sigma_i \le 2 \mu \\
	\hat{y}_i & 2\mu < \sigma_i < \theta \mu \;\text{and}\; h_2(\hat{y}_i) \le h_3(\theta \mu) \\
	\theta \mu & 2\mu < \sigma_i < \theta \mu \;\text{and}\; h_2(\hat{y}_i) > h_3(\theta \mu) \\
	\sigma_i & \sigma_i \ge \theta\mu
	\end{cases}
	\end{align*}
	depends on $\sigma_i$ and $\hat{y}_i = \frac{(\theta - 1)\sigma_i - \theta \mu}{\theta - 2}$.
\end{lemma}

\end{document}